\documentclass{amsart}
\usepackage{amssymb,amsmath,amsthm,latexsym,booktabs,todonotes, xcolor, comment, lineno}

\theoremstyle{definition}
\newtheorem{definition}{Definition}[section]

\theoremstyle{plain}
\newtheorem{lemma}[definition]{Lemma}

\newtheorem{theorem}[definition]{Theorem}

\usepackage{todonotes}

\usepackage[T1]{fontenc}
\usepackage{lmodern}
\usepackage{amsfonts,footnote,multicol,enumitem,multirow}
\usepackage[flushleft]{threeparttable}
\usepackage[font=small,labelfont=bf]{caption}
\usepackage{eucal}
\IfFileExists{esvect.sty}{\usepackage{esvect}}{\newcommand{\vv}[1]{\overrightarrow{##1}}}
\usepackage{graphicx,adjustbox}
\usepackage[noadjust]{cite}
\usepackage[hidelinks]{hyperref}
\numberwithin{equation}{section}
\theoremstyle{definition}
\newtheorem{ex}{Example}

\newtheorem{exb}{Example}

\newtheorem{exc}{Example}

\definecolor{orcidgreen}{HTML}{A6CE39}
\DeclareRobustCommand{\orcid}[1]{\href{https://orcid.org/#1}{\textsuperscript{\,\textcolor{orcidgreen}{\sffamily\bfseries iD}}}}
\newcommand{\Acknowledgements}[2]{\section*{Acknowledgements}#2\par}
\newcommand{\Competinginterests}[2]{\subsubsection*{Conflict of interest statement}#2\par}
\newcommand{\Availability}[2]{\subsubsection*{Data availability}#2\par}
\newcommand{\fitarray}[1]{\adjustbox{max width=\linewidth}{$\displaystyle #1$}}
\begin{document}
\title[The directed Hamilton--Waterloo problem]{The directed Hamilton--Waterloo problem for cycle lengths $3$ and $n$}
\author[U. Odaba\c{s}{\i}]{U\u{g}ur ODABA\c{S}I\orcid{0000-0001-6115-0342}}
\address{Department of Engineering Sciences, Istanbul University-Cerrahpa\c{s}a, Istanbul, 34320, T\"{u}rkiye}
\email{ugur.odabasi@iuc.edu.tr}
\author[D. Roberts]{Dan ROBERTS\orcid{0000-0001-8456-9533}}
\address{Illinois Wesleyan University, Bloomington, IL, 61701, USA}
\email{drobert1@iwu.edu}
\thanks{Dan Roberts is the corresponding author.}
\author[R.M. Low]{Richard M. LOW\orcid{0000-0003-2320-8296}}
\address{San Jose State University, San Jose, CA, 95192, USA}
\email{richard.low@sjsu.edu}
\date{}
\begin{abstract}
The directed Hamilton--Waterloo problem asks for a factorization of the complete symmetric digraph $K_v^*$ into directed $2$-factors, each of which is isomorphic to one of two given directed $2$-factors. In its uniform version, one seeks a decomposition of $K_v^*$ into $r$ factors consisting of directed $m$-cycles and $s$ factors consisting of directed $n$-cycles, where $r+s=v-1$. Let $m=3$, and let $n \geq 5$ be odd and not divisible by $3$. We show that, for every such $n$, the problem reduces to four decompositions of order at most $12n$: factorizations of $K^*_{3n}$ and $K^*_{6n}$ with every admissible number of factors of each kind, and two factorizations of the complete symmetric equipartite digraph with four parts of size $3n$. We construct these decompositions for $n=7$, and thereby give a complete solution for $(m,n)=(3,7)$: the obvious necessary conditions are also sufficient, with the single exception of $(v,r,s)=(6,5,0)$.
\end{abstract}
\maketitle

\section{Introduction}

A \emph{$k$-factor} of a graph $G$ is a spanning $k$-regular subgraph of $G$, and a \emph{$k$-factorization} of $G$ is a partition of its edge set into $k$-factors. A $2$-factor all of whose components are isomorphic to a graph $F$ is called an \emph{$F$-factor}. If the edge set of $G$ can be partitioned into $r_i$ factors consisting entirely of copies of $F_i$ for $1 \leq i \leq t$, we say that $G$ admits an $\{F_1^{r_1}, F_2^{r_2}, \ldots, F_t^{r_t}\}$-factorization; when every $F_i$ is a cycle $C_{n_i}$, this is a $\{C_{n_1}^{r_1}, C_{n_2}^{r_2}, \ldots, C_{n_t}^{r_t}\}$-factorization.

For a simple graph $G$, the \emph{symmetric digraph} $G^*$ has vertex set $V(G^*)=V(G)$ and arc set $A(G^*) = \bigcup_{\{y,z\} \in E(G)} \{(y,z),(z,y)\}$. In particular, the \emph{complete symmetric digraph} $K_v^*$ is the symmetric digraph of the complete graph $K_v$. The definitions above carry over to digraphs: a \emph{directed $2$-factor} of $G^*$ is a spanning subdigraph in which every vertex has in-degree and out-degree $1$, and a partition of the arc set of $G^*$ into directed $2$-factors is a \emph{directed $2$-factorization}. A directed $2$-factor all of whose components are directed $m$-cycles is a \emph{$\vv{C}_m$-factor}, and a partition of the arc set into $\vv{C}_m$-factors is a \emph{$\vv{C}_m$-factorization}. More generally, $G^*$ has a $\{\vv{C}_{n_1}^{\,r_1}, \ldots, \vv{C}_{n_t}^{\,r_t}\}$-factorization if its arc set can be partitioned into $r_i$ factors consisting entirely of directed $n_i$-cycles, $1 \leq i \leq t$.

The \emph{Oberwolfach problem} $\mathrm{OP}(m_1, \ldots, m_t)$ asks whether $v = m_1 + \cdots + m_t$ participants can be seated around $t$ round tables of sizes $m_1, \ldots, m_t$ on a number of nights so that every two participants are neighbors exactly once. In the \emph{directed} version $\mathrm{OP}^*(m_1, \ldots, m_t)$, every participant has to sit exactly once to the right of every other participant; equivalently, one asks for a decomposition of $K_v^*$ into $v-1$ directed $2$-factors, each consisting of directed cycles of lengths $m_1, \ldots, m_t$. The uniform directed problem $\mathrm{OP}^*(m^k)$, in which there are $k$ tables of size $m$, asks for a $\vv{C}_m$-factorization of $K_{km}^*$, and it is completely solved.

\begin{theorem}[\cite{Abel2002,Adamsun,Bennett1990,Bermond1979,Sajna2018,Sajna2014,Lacaze2023,Tillson1980}]\label{OPstar}
Let $m \geq 2$ and $k \geq 1$ be integers. Then $\mathrm{OP}^*(m^k)$ has a solution if and only if $(m,k) \notin \{(3,2),(4,1),(6,1)\}$.
\end{theorem}

The \emph{Hamilton--Waterloo problem} is the variant of the Oberwolfach problem in which the participants of a conference meet in two different venues, one in Hamilton and one in Waterloo, with different table arrangements. In graph-theoretic terms, it asks for a $2$-factorization of $K_v$ (or of $K_v - I$, where $I$ is a $1$-factor, when $v$ is even) in which every $2$-factor is isomorphic to one of two given $2$-factors. In the uniform case, $\mathrm{HWP}(v; m^r, n^s)$ asks for a factorization into $r$ $C_m$-factors and $s$ $C_n$-factors. Much work has been done for $m=3$; see, for instance, \cite{WCC} and the references therein. For $n=7$, Lei and Fu \cite{LeiFu} solved the case $v \equiv 21 \pmod{42}$ with three possible exceptions, all of order $21$, and Wang, Chen, and Cao \cite{WCC} settled these three cases together with the case $v \equiv 0 \pmod{42}$, thereby completing the solution of the Hamilton--Waterloo problem for triangle-factors and heptagon-factors.

In the \emph{directed Hamilton--Waterloo problem}, every participant must sit exactly once to the right of every other participant. Its uniform version $\mathrm{HWP}^*(v; m^r, n^s)$ asks for a $\{\vv{C}_m^{\,r}, \vv{C}_n^{\,s}\}$-factorization of $K_v^*$ with $r+s=v-1$, and such a factorization is denoted by $\mathrm{HW}^*(v; m^r, n^s)$. The following necessary conditions are immediate.

\begin{lemma}[\cite{FUS}]\label{NCs}
If an $\mathrm{HW}^*(v; m^{r}, n^{s})$ exists, then
\begin{enumerate}[label=\textup{(\roman*)}]
\item if $r>0$, then $v \equiv 0 \pmod m$;
\item if $s>0$, then $v \equiv 0 \pmod n$;
\item $r+s=v-1$.
\end{enumerate}
\end{lemma}

A directed almost resolvable version of the problem was considered by Dejter, Pike, and Rodger \cite{DPR}. For the directed Hamilton--Waterloo problem as stated above, Yetgin, Odaba\c{s}{\i}, and \"Ozkan \cite{FUS} settled, with a few possible exceptions, the cases $(m,n) \in \{(3,5),\allowbreak(3,15),\allowbreak(5,15),\allowbreak(4,6),\allowbreak(4,8),\allowbreak(4,12),\allowbreak(4,16),\allowbreak(6,12),\allowbreak(8,16)\}$. In \cite{FUS2}, the same authors considered factorizations of $K_v^*$ into factors of directed even cycles, namely $\mathrm{HW}^*(v; 2^r, m^s)$ and $\mathrm{HW}^*(v; m^r, (2m)^s)$ for even $m$. In \cite{YOK}, directed $4$-cycles and directed odd cycles were considered: the cases $(m,n) \in \{(3,4),(4,5),(4,7)\}$ were solved completely, and for every odd $n \geq 9$ the case $(m,n)=(4,n)$ was solved with a few possible exceptions.

In this paper, we consider the uniform directed Hamilton--Waterloo problem for $m=3$ and odd $n \geq 5$ with $3 \nmid n$. Our first result reduces the problem, for each such $n$, to four decompositions of small order; here $K^*_{(3n:4)}$ denotes the complete symmetric equipartite digraph with four parts of size $3n$ (see Section~\ref{sec:prelim}).

\begin{theorem}\label{reduction}
Let $n \geq 5$ be an odd integer with $3 \nmid n$, and suppose that
\begin{enumerate}[label=\textup{(P\arabic*)}]
\item an $\mathrm{HW}^*(3n; 3^{r}, n^{3n-1-r})$ exists for every integer $r$ with $0 \leq r \leq 3n-1$;
\item an $\mathrm{HW}^*(6n; 3^{r}, n^{6n-1-r})$ exists for every integer $r$ with $0 \leq r \leq 6n-1$;
\item $K^*_{(3n:4)}$ has a $\{\vv{C}_3^{\,a}, \vv{C}_n^{\,9n-a}\}$-factorization for $a = 3n$ and for $a = 6n$.
\end{enumerate}
Let $v \geq 2$ be an integer, and let $r$ and $s$ be nonnegative integers. Then an $\mathrm{HW}^*(v; 3^r, n^s)$ exists if and only if
\begin{enumerate}[label=\textup{(\roman*)}]
\item $v \equiv 0 \pmod 3$ when $r>0$,
\item $v \equiv 0 \pmod n$ when $s>0$, and
\item $r+s=v-1$,
\end{enumerate}
except when $(v,r,s)=(6,5,0)$.
\end{theorem}

For $n=7$, we establish (P1)--(P3) by explicit constructions, and obtain the following complete solution; it is the directed analogue of the result of \cite{LeiFu,WCC}.

\begin{theorem}\label{main}
Let $v \geq 2$ be an integer, and let $r$ and $s$ be nonnegative integers. Then an $\mathrm{HW}^*(v; 3^r, 7^s)$ exists if and only if
\begin{enumerate}[label=\textup{(\roman*)}]
\item $v \equiv 0 \pmod 3$ when $r>0$,
\item $v \equiv 0 \pmod 7$ when $s>0$, and
\item $r+s=v-1$,
\end{enumerate}
except when $(v,r,s)=(6,5,0)$.
\end{theorem}

The paper is organized as follows. Section~\ref{sec:prelim} collects the results that we use. In Section~\ref{sec:proof} we prove Theorem~\ref{reduction}, treating $v = 3nx$ with $x$ odd and $x$ even separately. In Section~\ref{sec:n7} we establish (P1)--(P3) for $n=7$ and prove Theorem~\ref{main}. The constructions for $K^*_{21}$ and $K^*_{42}$ are listed in Appendices~A and~B, and two factorizations of $K^*_{(21:4)}$ in Appendix~C.

\section{Preliminaries}\label{sec:prelim}

We use the following notation. For a digraph $D$ and a set $S \subseteq V(D)$, $D[S]$ denotes the subdigraph of $D$ induced by $S$. If $D_1, \ldots, D_k$ are pairwise arc-disjoint digraphs on the same vertex set, then $D_1 \oplus \cdots \oplus D_k$ denotes the digraph on this vertex set whose arc set is $A(D_1) \cup \cdots \cup A(D_k)$. For a positive integer $\alpha$ and a digraph $H$, $\alpha H$ denotes the vertex-disjoint union of $\alpha$ copies of $H$.

Undirected factorizations can be turned into directed ones as follows.

\begin{lemma}\label{doubling}
If $K_v$ has a $\{C_m^{a}, C_n^{b}\}$-factorization with $m,n \geq 3$, then an $\mathrm{HW}^*(v; m^{2a}, n^{2b})$ exists.
\end{lemma}

\begin{proof}
Let $F$ be a $C_k$-factor of $K_v$, $k \geq 3$. Orienting every cycle of $F$ in one of its two directions gives a $\vv{C}_k$-factor of $K_v^*$, and orienting every cycle in the opposite direction gives a second one; together they contain both arcs $(y,z)$ and $(z,y)$ for every edge $\{y,z\}$ of $F$, and no other arcs. Doing this for every factor of the given factorization, we obtain $2a$ $\vv{C}_m$-factors and $2b$ $\vv{C}_n$-factors whose arc sets partition the arc set of $K_v^*$.
\end{proof}

A \emph{Kirkman triple system} of order $v$, denoted by $\mathrm{KTS}(v)$, is a $C_3$-factorization of $K_v$; its $C_3$-factors are called \emph{parallel classes}, and the $3$-cycles of its parallel classes are called \emph{triangles}.

\begin{theorem}[\cite{RCW}]\label{KTS}
There exists a $\mathrm{KTS}(v)$ if and only if $v \equiv 3 \pmod 6$.
\end{theorem}

For integers $u, w \geq 2$, let $K_{(u:w)}$ denote the complete equipartite graph with $w$ parts of size $u$, and let $K_{(u:w)}^*$ be the corresponding complete symmetric equipartite digraph. Franceti\'c and \v{S}ajna \cite{Francetic2023} solved many cases of the directed uniform Oberwolfach problem for these digraphs.

\begin{theorem}[\cite{Francetic2023}]\label{dliu}
Let $u$, $w$, and $t$ be integers greater than $1$, and let $g=\gcd(w,t)$. Assume that one of the following conditions holds:
\begin{enumerate}[label=\textup{(\arabic*)}]
\item $u(w-1)$ is even;
\item $g \notin \{1,3\}$;
\item $g=1$, and $w \equiv 0 \pmod 4$ or $w \equiv 0 \pmod 6$;
\item $g=3$, and if $w=6$, then $u$ is divisible by a prime $p \leq 37$.
\end{enumerate}
Then $K_{(u:w)}^*$ has a $\vv{C}_t$-factorization if and only if $t \mid uw$, and $t$ is even in case $w=2$.
\end{theorem}

Every vertex of $K^*_{(b:3)}$ has out-degree $2b$, so every directed $2$-factorization of $K^*_{(b:3)}$ has exactly $2b$ factors. Similarly, every directed $2$-factorization of $K^*_{(3n:4)}$ has exactly $9n$ factors.

\begin{lemma}\label{K73}
Let $n \geq 5$ be odd with $3 \nmid n$, and let $b \in \{n,2n\}$. Then
\begin{enumerate}[label=\textup{(\alph*)}]
\item $K^*_{(b:3)}$ has a $\vv{C}_3$-factorization and a $\vv{C}_n$-factorization, each consisting of $2b$ factors;
\item $K^*_{(3n:4)}$ has a $\vv{C}_3$-factorization and a $\vv{C}_n$-factorization, each consisting of $9n$ factors.
\end{enumerate}
\end{lemma}

\begin{proof}
For (a), apply Theorem~\ref{dliu} with $u=b$ and $w=3$. Since $u(w-1)=2b$ is even, condition~(1) holds; moreover $t \in \{3,n\}$ divides $uw=3b$, and $w \neq 2$. For (b), apply Theorem~\ref{dliu} with $u=3n$ and $w=4$. Since $t \in \{3,n\}$ is odd, we have $\gcd(w,t)=1$, and $w \equiv 0 \pmod 4$, so condition~(3) holds; moreover $t$ divides $uw=12n$, and $w \neq 2$.
\end{proof}

We also need $C_3$-factorizations of complete equipartite graphs, that is, resolvable group divisible designs with block size $3$; see \cite[Theorem~1.1]{WCC} and the references therein.

\begin{theorem}[\cite{ASSW,AssafHartman,HoffmanSchellenberg,Liu2000,Liu2003,Piotrowski,Rees93}]\label{RGDD}
Let $u, w \geq 2$. Then $K_{(u:w)}$ has a $C_3$-factorization if and only if $w \geq 3$, $u(w-1) \equiv 0 \pmod 2$, $uw \equiv 0 \pmod 3$ and $(u,w) \notin \{(2,3),(2,6),(6,3)\}$.
\end{theorem}

\section{The reduction}\label{sec:proof}

Throughout this section, $n \geq 5$ is an odd integer with $3 \nmid n$, and we assume that (P1)--(P3) hold.

\begin{lemma}\label{K84}
Let $r$ be an integer with $0 \leq r \leq 12n-1$. Then an $\mathrm{HW}^*(12n; 3^{r}, n^{12n-1-r})$ exists.
\end{lemma}

\begin{proof}
Partition the vertex set of $K^*_{12n}$ into four sets $W_1, \ldots, W_4$ of size $3n$. The arcs of $K^*_{12n}$ with both ends in the same set $W_j$ form a spanning subdigraph isomorphic to $4K^*_{3n}$, and the arcs joining two different sets $W_j$ form a spanning subdigraph $D \cong K^*_{(3n:4)}$; thus $K^*_{12n} \cong 4K^*_{3n} \oplus K^*_{(3n:4)}$. Put $a = 3n \lfloor r/(3n) \rfloor$ and $r_0 = r - a$; since $r \leq 12n-1$, we have $a \in \{0,3n,6n,9n\}$ and $0 \leq r_0 \leq 3n-1$. On each $K^*_{12n}[W_j] \cong K^*_{3n}$ take an $\mathrm{HW}^*(3n; 3^{r_0}, n^{3n-1-r_0})$, which exists by (P1); since $W_1, \ldots, W_4$ are pairwise disjoint and cover all vertices, taking the unions of corresponding factors over $j$ (matching $\vv{C}_3$-factors with $\vv{C}_3$-factors and $\vv{C}_n$-factors with $\vv{C}_n$-factors) gives $r_0$ $\vv{C}_3$-factors and $3n-1-r_0$ $\vv{C}_n$-factors of $K^*_{12n}$. By Lemma~\ref{K73}(b) if $a \in \{0,9n\}$, and by (P3) if $a \in \{3n,6n\}$, $D$ has a $\{\vv{C}_3^{\,a}, \vv{C}_n^{\,9n-a}\}$-factorization. Together these form an $\mathrm{HW}^*(12n; 3^{r_0+a}, n^{12n-1-r_0-a})$, and $r_0 + a = r$.
\end{proof}

All other constructions in this section are instances of the following lemma, in which a resolvable design with block size $3$ is blown up and its groups are filled with small directed Hamilton--Waterloo designs.

\begin{lemma}\label{lifting}
Let $u \geq 2$, $w \geq 3$, and $b \in \{n,2n\}$ be integers. Suppose that $K_{(u:w)}$ has a $C_3$-factorization and that an $\mathrm{HW}^*(ub; 3^{r_0}, n^{ub-1-r_0})$ exists for every integer $r_0$ with $0 \leq r_0 \leq ub-1$. Then an $\mathrm{HW}^*(ubw; 3^{r}, n^{s})$ exists for all nonnegative integers $r$ and $s$ with $r+s = ubw-1$.
\end{lemma}

\begin{proof}
Let $X_1, \ldots, X_w$ be the parts of $K_{(u:w)}$, and let $\mathcal{T}_1, \ldots, \mathcal{T}_q$ be the parallel classes of the given $C_3$-factorization; since every vertex has degree $u(w-1)$ and each parallel class contains exactly two edges at every vertex, we have $q = u(w-1)/2$.

Blow up every vertex $x$ of $K_{(u:w)}$ into a set $P_x$ of $b$ vertices, and put $W_h = \bigcup_{x \in X_h} P_x$, so that $|W_h| = ub$. For a triangle $T = \{x,y,z\}$ of the $C_3$-factorization, let $D_T \cong K^*_{(b:3)}$ be the subdigraph of $K^*_{ubw}$ formed by the arcs joining two different sets among $P_x$, $P_y$ and $P_z$. Every arc of $K^*_{ubw}$ either joins two vertices of the same set $W_h$, or joins vertices of $P_x$ and $P_y$ for an edge $xy$ of $K_{(u:w)}$, and then it lies in $D_T$ for the unique triangle $T$ containing $xy$. Hence
\[
K^*_{ubw} = H_0 \oplus H_1 \oplus \cdots \oplus H_q
\;\cong\; wK^*_{ub} \oplus
\underbrace{\tfrac{uw}{3} K^*_{(b:3)} \oplus \cdots \oplus
\tfrac{uw}{3} K^*_{(b:3)}}_{q},
\]
where $H_0$ consists of the arcs within the sets $W_h$, and $H_j$, for $1 \leq j \leq q$, of the arcs of the digraphs $D_T$ with $T \in \mathcal{T}_j$; since $\mathcal{T}_j$ is a parallel class, these digraphs are pairwise vertex-disjoint and span all vertices.

Let $0 \leq r_0 \leq ub-1$ and $0 \leq k \leq q$. Take an $\mathrm{HW}^*(ub; 3^{r_0}, n^{ub-1-r_0})$ on every component of $H_0$, and, by Lemma~\ref{K73}(a), a $\vv{C}_3$-factorization of every component of $H_j$ for $1 \leq j \leq k$ and a $\vv{C}_n$-factorization of every component of $H_j$ for $k < j \leq q$, each with $2b$ factors. Uniting factors of the same kind across the components of each $H_j$, $0 \leq j \leq q$, we obtain a directed $2$-factorization of $K^*_{ubw}$ with $r = r_0 + 2bk$ $\vv{C}_3$-factors and, since $(ub-1) + 2bq = ubw-1$, with $ubw-1-r$ $\vv{C}_n$-factors. For a fixed $k$, these values of $r$ form the interval $[2bk,\, 2bk+ub-1]$. Since $u \geq 2$, we have $2b \leq ub$, so consecutive intervals overlap or are adjacent, and their union is $[0,\, 2bq+ub-1] = [0,\, ubw-1]$.
\end{proof}

\subsection{When \texorpdfstring{$x$}{x} is odd}\label{sec:odd}

\begin{theorem}\label{xodd}
Let $x$ be an odd positive integer and let $v = 3nx$. Then an $\mathrm{HW}^*(v; 3^{r}, n^{s})$ exists for all nonnegative integers $r$ and $s$ with $r + s = v - 1$.
\end{theorem}

\begin{proof}
The case $x=1$ is (P1), so let $x \geq 3$. Since $x$ is odd, $3x \equiv 3 \pmod 6$, and by Theorem~\ref{KTS} there is a $\mathrm{KTS}(3x)$. Remove one of its parallel classes and regard its $x$ triangles as the parts of a complete equipartite graph $K_{(3:x)}$. Every edge of $K_{3x}$ lies in exactly one triangle of the $\mathrm{KTS}(3x)$. The edges within the parts lie in the triangles of the removed class, so no triangle of the remaining $(3x-3)/2$ parallel classes contains such an edge, and every edge of $K_{(3:x)}$ lies in exactly one of these triangles. Since each of these classes spans all vertices, they form a $C_3$-factorization of $K_{(3:x)}$. By (P1), the result follows from Lemma~\ref{lifting} with $u=3$, $w=x$, and $b=n$.
\end{proof}

\subsection{When \texorpdfstring{$x$}{x} is even}\label{sec:even}

When $x$ is even, $K_{3x}$ has odd degree, so there is no $\mathrm{KTS}(3x)$. The natural substitute, a $C_3$-factorization of $K_{3x}-I$ for a $1$-factor $I$, does not help: blowing up an edge of $I$ would produce a copy of $K^*_{(n:2)}$, which is bipartite and therefore has no directed cycles of odd length. Instead, we apply Lemma~\ref{lifting} with $ub = 6n$.

\begin{theorem}\label{xeven}
Let $x \neq 4$ be an even positive integer and let $v = 3nx$. Then an $\mathrm{HW}^*(v; 3^{r}, n^{s})$ exists for all nonnegative integers $r$ and $s$ with $r + s = v - 1$.
\end{theorem}

\begin{proof}
Write $x = 2y$. The case $y=1$ is (P2), so let $y \geq 3$. If $y \geq 4$, then $K_{(6:y)}$ has a $C_3$-factorization by Theorem~\ref{RGDD}, and we apply Lemma~\ref{lifting} with $u = 6$, $w = y$ and $b = n$. If $y = 3$, then $K_{(3:3)}$ has a $C_3$-factorization by Theorem~\ref{RGDD}, and we apply Lemma~\ref{lifting} with $u = 3$, $w = 3$ and $b = 2n$. In both cases $ub = 6n$, and the required $\mathrm{HW}^*(6n; 3^{r_0}, n^{6n-1-r_0})$ exist by (P2).
\end{proof}

Note that the order $v=12n$, which corresponds to $x=4$ and $y=2$, is not covered by Theorem~\ref{xeven}: Lemma~\ref{lifting} requires $w \geq 3$, but $ub = 6n$ and $v = 12n$ would force $w = 2$, and $K_{(u:2)}$ is bipartite and has no $C_3$-factorization. This order is settled by Lemma~\ref{K84} instead.

\subsection{Proof of Theorem~\ref{reduction}}

The necessity of the conditions follows from Lemma~\ref{NCs}. For sufficiency, first let $r=0$ or $s=0$. Then the problem is $\mathrm{OP}^*(n^{v/n})$ or $\mathrm{OP}^*(3^{v/3})$, which by Theorem~\ref{OPstar} has a solution unless $(3, v/3)=(3,2)$, that is, unless $(v,r,s)=(6,5,0)$. So let $r,s \geq 1$. By Lemma~\ref{NCs}, and since $\gcd(3,n)=1$, we have $v \equiv 0 \pmod{3n}$; write $v = 3nx$ with $x \geq 1$. If $x$ is odd, the result follows from Theorem~\ref{xodd}; if $x$ is even and $x \neq 4$, from Theorem~\ref{xeven}; and if $x=4$, that is, $v=12n$, from Lemma~\ref{K84}. \qed

\section{The case \texorpdfstring{$n=7$}{n=7}}\label{sec:n7}

In this section we establish (P1)--(P3) for $n=7$, that is, for the orders $3n = 21$ and $6n = 42$ and for the digraph $K^*_{(21:4)}$.

The undirected Hamilton--Waterloo problem for triangle-factors and heptagon-factors has been solved completely.

\begin{theorem}[\cite{LeiFu,WCC}]\label{HWP37}
Let $v \equiv 0 \pmod{21}$, and let $a$ and $b$ be nonnegative integers with $a+b = \lfloor \frac{v-1}{2} \rfloor$. Then $K_v$ if $v$ is odd, and $K_v$ minus a $1$-factor if $v$ is even, has a $\{C_3^{a}, C_7^{b}\}$-factorization.
\end{theorem}

By Lemma~\ref{doubling}, an undirected $\{C_3^{a}, C_7^{b}\}$-factorization of $K_v$ only yields directed factorizations in which both $r=2a$ and $s=2b$ are even. For $v=21$ we have $r+s=20$, so $r$ and $s$ have the same parity, and the cases with $r$ odd require new constructions.

\begin{lemma}\label{K21}
For every integer $r$ with $0 \leq r \leq 20$, there exists an $\mathrm{HW}^*(21; 3^{r}, 7^{20-r})$.
\end{lemma}

\begin{proof}
If $r$ is even, then by Theorem~\ref{HWP37} with $v=21$, $a=r/2$ and $b=(20-r)/2$, the complete graph $K_{21}$ has a $\{C_3^{r/2}, C_7^{(20-r)/2}\}$-factorization, and Lemma~\ref{doubling} yields an $\mathrm{HW}^*(21; 3^{r}, 7^{20-r})$. If $r$ is odd, that is, $r \in \{1,3,\ldots,19\}$, see Examples~\ref{hw21r1s19}--\ref{hw21r19s1} in Appendix~A.
\end{proof}

The constructions in Appendix~A are given in the usual way, by base cycles on $\mathbb{Z}_{21}$ developed by a cyclic automorphism group of order $7$ or $3$.

For $v=42$, Lemma~\ref{doubling} does not apply at all: $K_{42}$ has odd degree and therefore no $2$-factorization, and Theorem~\ref{HWP37} only gives factorizations of $K_{42}$ minus a $1$-factor, whose doubling does not cover the arcs joining the endpoints of the edges of this $1$-factor. Hence all cases with $1 \leq r \leq 40$ require new constructions.

\begin{lemma}\label{K42}
For every integer $r$ with $0 \leq r \leq 41$, there exists an $\mathrm{HW}^*(42; 3^{r}, 7^{41-r})$.
\end{lemma}

\begin{proof}
For $r=0$ and $r=41$, the required decompositions are a $\vv{C}_7$-factorization and a $\vv{C}_3$-factorization of $K^*_{42}$, respectively; they exist by Theorem~\ref{OPstar}, since $(7,6),(3,14) \notin \{(3,2),(4,1),(6,1)\}$. For $1 \leq r \leq 40$, see Examples~\ref{hw42r1s40}--\ref{hw42r40s1} in Appendix~B, where all constructions are developed by the cyclic group $C_{21}$ acting on $\mathbb{Z}_{42}$ by $x \mapsto x+2$.
\end{proof}

\begin{lemma}\label{K214}
For $a \in \{21, 42\}$, $K^*_{(21:4)}$ has a $\{\vv{C}_3^{\,a}, \vv{C}_7^{\,63-a}\}$-factorization.
\end{lemma}

\begin{proof}
See Examples~\ref{k214a21} and~\ref{k214a42} in Appendix~C.
\end{proof}

\begin{proof}[Proof of Theorem~\ref{main}]
For $n=7$, conditions (P1), (P2) and (P3) of Theorem~\ref{reduction} hold by Lemmas~\ref{K21}, \ref{K42}, and~\ref{K214}, respectively. Hence Theorem~\ref{main} follows from Theorem~\ref{reduction}.
\end{proof}

\section{Concluding remarks}

Theorem~\ref{reduction} reduces the uniform directed Hamilton--Waterloo problem for cycle lengths $3$ and $n$, where $n \geq 5$ is odd and $3 \nmid n$, to the four decompositions in (P1)--(P3), all of order at most $12n$. For $n=7$ we have constructed these decompositions, and the necessary conditions of Lemma~\ref{NCs} are therefore sufficient for the existence of an $\mathrm{HW}^*(v; 3^r, 7^s)$, the only exception being $(v,r,s)=(6,5,0)$, which comes from the nonexistence of a solution to $\mathrm{OP}^*(3^2)$. The proof combines Kirkman triple systems, the results of Franceti\'c and \v{S}ajna \cite{Francetic2023} on complete symmetric equipartite digraphs, resolvable group divisible designs, the solution of the undirected problem \cite{LeiFu,WCC}, and explicit decompositions of $K^*_{21}$, $K^*_{42}$ and $K^*_{(21:4)}$ obtained from cyclic group actions.

For other values of $n$, establishing (P1)--(P3) remains open. For $n=5$, the case of odd $v$ was considered in \cite{FUS}; by Theorem~\ref{reduction}, a solution for all $v$ would follow from (P1)--(P3) for $n=5$, that is, from suitable decompositions of $K^*_{15}$, $K^*_{30}$ and $K^*_{(15:4)}$.

\Acknowledgements{}{During manuscript preparation, the authors used ChatGPT 6.0 Astra for improving expositional clarity and general proofreading.  It was also used to generate Python code for constructing tables in the appendices. The final manuscript was reviewed and verified by the authors. The tool was not used for research design, data generation or analysis, interpretation, analytical code, mathematical proofs, or scientific conclusions. The authors take full responsibility for the final manuscript.}

\Competinginterests{}{The authors declare that they have no known competing financial interests or personal relationships that could have appeared to
	influence the work reported in this paper.}

\Availability{}{All decompositions used in the proofs are listed explicitly in Appendices~A, B, and~C.}

\appendix
\section*{Appendix}
\addcontentsline{toc}{section}{Appendix}

In what follows, $(a_0, a_1, \ldots, a_{k-1})$ with $k \in \{3,7\}$ denotes the directed $k$-cycle with arcs $(a_i, a_{i+1})$ for $0 \leq i \leq k-1$, where the subscripts are taken modulo $k$. In addition, in the tables and in the descriptions of the automorphism groups, $C_n$ denotes the cyclic group of order $n$ (not a cycle), and $1$ denotes the identity subgroup. In each construction, the base cycles given in the table are developed by the indicated automorphism group. The base cycles in each row of a table, when developed by the subgroup given in the last column of the corresponding row, yield a cycle factor; and the action of each coset of this subgroup on these base cycles yields a different cycle factor.

\subsection*{A. The order 21}

\begin{ex}\label{hw21r1s19} $\mathrm{HW}^*(21;3^{1},7^{19})$ on $\mathbb{Z}_{21}$ developed by $C_{7}=\left<\alpha\right>$ where $\alpha~:~x\rightarrow x+3\pmod{21}$. There is 1 orbit of length 7 on 3-cycles, 6 orbits of length 7 on 7-cycles, and 15 fixed 7-cycles.
\[\fitarray{\begin{array}{|l l l|c|}
\hline
(0,7,14) & & & C_7 \\
\hline
(0,17,4,6,11,1,5) & (2,19,8,18,7,15,10) & (3,16,14,12,13,9,20) & \multirow{2}{*}{$1$} \\
(0,18,1,16,17,5,7) & (2,3,11,4,20,19,9) & (6,8,13,12,10,15,14) & \\
\hline
(0,3,6,9,12,15,18) & (1,4,7,10,13,16,19) & (2,5,8,11,14,17,20) & \\
(0,6,12,18,3,9,15) & (1,7,13,19,4,10,16) & (2,8,14,20,5,11,17) & \\
(0,9,18,6,15,3,12) & (1,10,19,7,16,4,13) & (2,14,5,17,8,20,11) & C_7 \\
(0,12,3,15,6,18,9) & (1,13,4,16,7,19,10) & (2,17,11,5,20,14,8) & \\
(0,15,9,3,18,12,6) & (1,19,16,13,10,7,4) & (2,20,17,14,11,8,5) & \\
\hline
\end{array}}\]
\end{ex}

\begin{ex}\label{hw21r3s17} $\mathrm{HW}^*(21;3^{3},7^{17})$ on $\mathbb{Z}_{21}$ developed by $C_{7}=\left<\alpha\right>$ where $\alpha~:~x\rightarrow x+3\pmod{21}$. There are 3 orbits of length 7 on 3-cycles, 6 orbits of length 7 on 7-cycles, and 9 fixed 7-cycles.
\[\fitarray{\begin{array}{|l l l |c|}
\hline
(0,2,19) &  &  & \multirow{3}{*}{$C_{7}$} \\
(0,13,2) &  &  &  \\
(0,17,16) &  &  &  \\
\hline
(0,20,2,16,1,15,9) & (3,10,18,11,19,13,4) & (5,7,8,14,6,12,17) & \multirow{2}{*}{$1$} \\
(0,16,8,13,11,12,10) & (1,5,9,19,15,2,18) & (3,4,20,6,17,7,14) &  \\
\hline
(0,9,18,6,15,3,12) & (1,10,19,7,16,4,13) & (2,17,11,5,20,14,8) & \multirow{3}{*}{$C_{7}$} \\
(0,3,6,9,12,15,18) & (1,4,7,10,13,16,19) & (2,14,5,17,8,20,11) &  \\
(0,18,15,12,9,6,3) & (1,19,16,13,10,7,4) & (2,20,17,14,11,8,5) &  \\
\hline
\end{array}}\]
\end{ex}

\begin{ex}\label{hw21r5s15} $\mathrm{HW}^*(21;3^{5},7^{15})$ on $\mathbb{Z}_{21}$ developed by $C_{7}=\left<\alpha\right>$ where $\alpha~:~x\rightarrow x+3\pmod{21}$. There are 5 orbits of length 7 on 3-cycles, 6 orbits of length 7 on 7-cycles, and 3 fixed 7-cycles.
\[\fitarray{\begin{array}{|l l l |c|}
\hline
(0,2,19) &  &  & \multirow{5}{*}{$C_{7}$} \\
(0,13,2) &  &  &  \\
(0,17,16) &  &  &  \\
(0,10,17) &  &  &  \\
(0,1,20) &  &  &  \\
\hline
(0,4,15,12,11,6,13) & (1,18,3,17,2,16,10) & (5,8,19,7,20,9,14) & \multirow{2}{*}{$1$} \\
(0,19,16,1,17,14,6) & (2,9,18,8,10,3,15) & (4,5,13,12,7,11,20) &  \\
\hline
(0,3,6,9,12,15,18) & (1,4,7,10,13,16,19) & (2,17,11,5,20,14,8) & C_{7} \\
\hline
\end{array}}\]
\end{ex}

\begin{ex}\label{hw21r7s13} $\mathrm{HW}^*(21;3^{7},7^{13})$ on $\mathbb{Z}_{21}$ developed by $C_{7}=\left<\alpha\right>$ where $\alpha~:~x\rightarrow x+3\pmod{21}$. There are 7 orbits of length 7 on 3-cycles, 3 orbits of length 7 on 7-cycles, and 18 fixed 7-cycles.
\[\fitarray{\begin{array}{|l l l |c|}
\hline
(0,11,1) & (2,7,3) & (4,17,18) & \multirow{3}{*}{$1$} \\
(5,15,16) & (6,20,13) & (8,10,12) &  \\
(9,19,14) &  &  &  \\
\hline
(0,16,20,7,8,4,2) & (1,6,19,5,12,10,18) & (3,11,15,17,9,14,13) & 1 \\
\hline
(0,3,6,9,12,15,18) & (1,13,4,16,7,19,10) & (2,17,11,5,20,14,8) & \multirow{6}{*}{$C_{7}$} \\
(0,9,18,6,15,3,12) & (1,16,10,4,19,13,7) & (2,11,20,8,17,5,14) &  \\
(0,18,15,12,9,6,3) & (1,7,13,19,4,10,16) & (2,8,14,20,5,11,17) &  \\
(0,15,9,3,18,12,6) & (1,10,19,7,16,4,13) & (2,5,8,11,14,17,20) &  \\
(0,12,3,15,6,18,9) & (1,4,7,10,13,16,19) & (2,20,17,14,11,8,5) &  \\
(0,6,12,18,3,9,15) & (1,19,16,13,10,7,4) & (2,14,5,17,8,20,11) &  \\
\hline
\end{array}}\]
\end{ex}

\begin{ex}\label{hw21r9s11} $\mathrm{HW}^*(21;3^{9},7^{11})$ on $\mathbb{Z}_{21}$ developed by $C_{3}=\left<\alpha\right>$ where $\alpha~:~x\rightarrow x+7\pmod{21}$. There are 21 orbits of length 3 on 3-cycles, and 11 orbits of length 3 on 7-cycles.
\[\fitarray{\begin{array}{|l l l |c|}
\hline
(0,15,4) & (1,18,2) & (3,7,5) & \multirow{3}{*}{$1$} \\
(6,14,13) & (8,9,17) & (10,12,20) &  \\
(11,19,16) &  &  &  \\
\hline
(0,4,18) & (1,10,17) & (2,5,9) & \multirow{3}{*}{$1$} \\
(3,13,14) & (6,8,16) & (7,20,15) &  \\
(11,12,19) &  &  &  \\
\hline
(0,9,16) & (1,17,12) & (2,3,6) & \multirow{3}{*}{$1$} \\
(4,13,10) & (5,7,14) & (8,19,15) &  \\
(11,18,20) &  &  &  \\
\hline
(0,1,16,20,19,18,10) &  &  & \multirow{2}{*}{$C_{3}$} \\
(0,6,9,11,5,17,15) &  &  &  \\
\hline
(0,2,12,6,13,11,9) & (1,4,8,7,3,17,14) & (5,16,10,18,15,19,20) & \multirow{3}{*}{$1$} \\
(0,11,1,20,4,3,12) & (2,19,14,17,18,9,7) & (5,8,15,6,10,16,13) &  \\
(0,18,13,1,3,2,11) & (4,10,6,16,8,5,19) & (7,12,17,9,15,20,14) &  \\
\hline
\end{array}}\]
\end{ex}

\begin{ex}\label{hw21r11s9} $\mathrm{HW}^*(21;3^{11},7^{9})$ on $\mathbb{Z}_{21}$ developed by $C_{7}=\left<\alpha\right>$ where $\alpha~:~x\rightarrow x+3\pmod{21}$. There are 11 orbits of length 7 on 3-cycles, 3 orbits of length 7 on 7-cycles, and 6 fixed 7-cycles.
\[\fitarray{\begin{array}{|l l l |c|}
\hline
(0,2,19) &  &  & \multirow{4}{*}{$C_{7}$} \\
(0,13,2) &  &  &  \\
(0,17,16) &  &  &  \\
(0,10,17) &  &  &  \\
\hline
(0,4,18) & (1,16,3) & (2,13,10) & \multirow{3}{*}{$1$} \\
(5,20,6) & (7,11,14) & (8,9,17) &  \\
(12,19,15) &  &  &  \\
\hline
(0,5,14,16,1,17,9) & (2,12,18,8,13,3,4) & (6,20,15,10,11,19,7) & 1 \\
\hline
(0,15,9,3,18,12,6) & (1,13,4,16,7,19,10) & (2,20,17,14,11,8,5) & \multirow{2}{*}{$C_{7}$} \\
(0,9,18,6,15,3,12) & (1,4,7,10,13,16,19) & (2,8,14,20,5,11,17) &  \\
\hline
\end{array}}\]
\end{ex}

\begin{ex}\label{hw21r13s7} $\mathrm{HW}^*(21;3^{13},7^{7})$ on $\mathbb{Z}_{21}$ developed by $C_{7}=\left<\alpha\right>$ where $\alpha~:~x\rightarrow x+3\pmod{21}$. There are 13 orbits of length 7 on 3-cycles, and 3 orbits of length 7 on 7-cycles.
\[\fitarray{\begin{array}{|l l l |c|}
\hline
(0,2,19) &  &  & \multirow{6}{*}{$C_{7}$} \\
(0,13,2) &  &  &  \\
(0,17,16) &  &  &  \\
(0,10,17) &  &  &  \\
(0,1,20) &  &  &  \\
(0,16,11) &  &  &  \\
\hline
(0,12,5) & (1,9,15) & (2,17,3) & \multirow{3}{*}{$1$} \\
(4,18,16) & (6,10,7) & (8,19,11) &  \\
(13,14,20) &  &  &  \\
\hline
(0,18,5,17,20,8,10) & (1,4,19,2,15,3,14) & (6,11,16,7,13,9,12) & 1 \\
\hline
\end{array}}\]
\end{ex}

\begin{ex}\label{hw21r15s5} $\mathrm{HW}^*(21;3^{15},7^{5})$ on $\mathbb{Z}_{21}$ developed by $C_{7}=\left<\alpha\right>$ where $\alpha~:~x\rightarrow x+3\pmod{21}$. There are 15 orbits of length 7 on 3-cycles, and 15 fixed 7-cycles.
\[\fitarray{\begin{array}{|l l l |c|}
\hline
(0,5,10) &  &  & \multirow{8}{*}{$C_{7}$} \\
(0,19,2) &  &  &  \\
(0,20,13) &  &  &  \\
(0,17,19) &  &  &  \\
(0,11,7) &  &  &  \\
(0,10,17) &  &  &  \\
(0,4,20) &  &  &  \\
(0,16,5) &  &  &  \\
\hline
(0,2,1) & (3,16,7) & (4,9,17) & \multirow{3}{*}{$1$} \\
(5,15,12) & (6,13,14) & (8,20,10) &  \\
(11,18,19) &  &  &  \\
\hline
(0,15,9,3,18,12,6) & (1,19,16,13,10,7,4) & (2,5,8,11,14,17,20) & \multirow{5}{*}{$C_{7}$} \\
(0,12,3,15,6,18,9) & (1,10,19,7,16,4,13) & (2,8,14,20,5,11,17) &  \\
(0,6,12,18,3,9,15) & (1,7,13,19,4,10,16) & (2,17,11,5,20,14,8) &  \\
(0,3,6,9,12,15,18) & (1,16,10,4,19,13,7) & (2,20,17,14,11,8,5) &  \\
(0,9,18,6,15,3,12) & (1,4,7,10,13,16,19) & (2,11,20,8,17,5,14) &  \\
\hline
\end{array}}\]
\end{ex}

\begin{ex}\label{hw21r17s3} $\mathrm{HW}^*(21;3^{17},7^{3})$ on $\mathbb{Z}_{21}$ developed by $C_{7}=\left<\alpha\right>$ where $\alpha~:~x\rightarrow x+3\pmod{21}$. There are 17 orbits of length 7 on 3-cycles, and 9 fixed 7-cycles.
\[\fitarray{\begin{array}{|l l l |c|}
\hline
(0,10,14) &  &  & \multirow{10}{*}{$C_{7}$} \\
(0,2,4) &  &  &  \\
(0,4,20) &  &  &  \\
(0,7,8) &  &  &  \\
(0,16,5) &  &  &  \\
(0,19,17) &  &  &  \\
(0,14,19) &  &  &  \\
(0,5,1) &  &  &  \\
(0,20,13) &  &  &  \\
(0,8,7) &  &  &  \\
\hline
(0,17,6) & (1,19,11) & (2,10,16) & \multirow{3}{*}{$1$} \\
(3,4,13) & (5,8,14) & (7,12,15) &  \\
(9,20,18) &  &  &  \\
\hline
(0,9,18,6,15,3,12) & (1,13,4,16,7,19,10) & (2,17,11,5,20,14,8) & \multirow{3}{*}{$C_{7}$} \\
(0,18,15,12,9,6,3) & (1,16,10,4,19,13,7) & (2,20,17,14,11,8,5) &  \\
(0,6,12,18,3,9,15) & (1,4,7,10,13,16,19) & (2,11,20,8,17,5,14) &  \\
\hline
\end{array}}\]
\end{ex}

\begin{ex}\label{hw21r19s1} $\mathrm{HW}^*(21;3^{19},7^{1})$ on $\mathbb{Z}_{21}$ developed by $C_{7}=\left<\alpha\right>$ where $\alpha~:~x\rightarrow x+3\pmod{21}$. There are 19 orbits of length 7 on 3-cycles, and 3 fixed 7-cycles.
\[\fitarray{\begin{array}{|l l l |c|}
\hline
(0,2,19) &  &  & \multirow{5}{*}{$C_{7}$} \\
(0,13,2) &  &  &  \\
(0,17,16) &  &  &  \\
(0,10,11) &  &  &  \\
(0,4,17) &  &  &  \\
\hline
(0,6,13) & (1,18,19) & (2,9,8) & \multirow{3}{*}{$1$} \\
(3,14,11) & (4,15,20) & (5,16,7) &  \\
(10,17,12) &  &  &  \\
\hline
(0,9,6) & (1,15,8) & (2,5,17) & \multirow{3}{*}{$1$} \\
(3,11,12) & (4,20,7) & (10,14,16) &  \\
(13,19,18) &  &  &  \\
\hline
(0,3,6,9,12,15,18) & (1,10,19,7,16,4,13) & (2,11,20,8,17,5,14) & C_{7} \\
\hline
\end{array}}\]
\end{ex}

\subsection*{B. The order 42}

All examples in this part are on $\mathbb{Z}_{42}$ and are developed by the cyclic group $C_{21}=\left<\alpha\right>$, where $\alpha~:~x\rightarrow x+2\pmod{42}$; the subgroups $\left<\alpha^{3}\right>$ and $\left<\alpha^{7}\right>$ that appear in the tables have orders $7$ and $3$, respectively.

\begin{exb}\label{hw42r1s40} $\mathrm{HW}^*(42;3^{1},7^{40})$. There are 2 orbits of length 7 on 3-cycles, and 10 orbits of length 21 and 10 orbits of length 3 on 7-cycles.
\setlength{\arraycolsep}{3pt}
\[\fitarray{\begin{array}{|l l l |c|}
\hline
(0,28,14) & (1,15,29) &  & C_{21} \\
\hline
(0,30,18,6,36,24,12) & (1,25,7,31,13,37,19) &  & \multirow{5}{*}{$C_{21}$} \\
(0,36,30,24,18,12,6) & (1,31,19,7,37,25,13) &  &  \\
(0,12,24,36,6,18,30) & (1,7,13,19,25,31,37) &  &  \\
(0,18,36,12,30,6,24) & (1,37,31,25,19,13,7) &  &  \\
(0,6,12,18,24,30,36) & (1,13,25,37,7,19,31) &  &  \\
\hline
(0,15,35,27,37,12,10) & (2,6,11,4,17,33,36) &  & \multirow{2}{*}{$\left<\alpha^{7}\right>$} \\
(0,20,12,8,11,38,13) & (1,3,2,4,23,21,5) &  &  \\
\hline
(0,29,40,5,38,22,33) & (1,9,28,19,16,25,11) & (2,37,18,34,17,12,39) & \multirow{2}{*}{$1$} \\
(3,21,4,14,36,35,24) & (6,20,15,30,31,10,41) & (7,8,32,13,26,23,27) &  \\
\hline
\end{array}}\]
\end{exb}

\begin{exb}\label{hw42r2s39} $\mathrm{HW}^*(42;3^{2},7^{39})$. There are 4 orbits of length 7 on 3-cycles, and 10 orbits of length 21 and 8 orbits of length 3 on 7-cycles.
\[\fitarray{\begin{array}{|l l |c|}
\hline
(0,14,28) & (1,15,29) & \multirow{2}{*}{$C_{21}$} \\
(0,28,14) & (1,29,15) &  \\
\hline
(0,36,30,24,18,12,6) & (1,25,7,31,13,37,19) & \multirow{4}{*}{$C_{21}$} \\
(0,24,6,30,12,36,18) & (1,13,25,37,7,19,31) &  \\
(0,30,18,6,36,24,12) & (1,7,13,19,25,31,37) &  \\
(0,6,12,18,24,30,36) & (1,31,19,7,37,25,13) &  \\
\hline
(0,37,34,19,15,24,13) & (2,35,26,3,4,39,22) & \multirow{5}{*}{$\left<\alpha^{7}\right>$} \\
(0,10,9,7,22,29,34) & (2,5,4,25,17,12,41) &  \\
(0,16,6,31,10,33,23) & (1,8,12,4,21,41,25) &  \\
(0,1,9,3,5,40,11) & (2,22,24,21,34,18,27) &  \\
(0,12,17,2,20,35,4) & (1,24,22,33,13,23,39) &  \\
\hline
\end{array}}\]
\end{exb}

\begin{exb}\label{hw42r3s38} $\mathrm{HW}^*(42;3^{3},7^{38})$. There are 2 orbits of length 21 on 3-cycles, and 10 orbits of length 21 and 6 orbits of length 3 on 7-cycles.
\setlength{\arraycolsep}{3pt}
\[\fitarray{\begin{array}{|l l l |c|}
\hline
(0,20,28) & (1,15,17) &  & \left<\alpha^{3}\right> \\
\hline
(0,6,12,18,24,30,36) & (1,7,13,19,25,31,37) & (2,26,8,32,14,38,20) & \multirow{2}{*}{$\left<\alpha^{3}\right>$} \\
(3,27,9,33,15,39,21) & (4,16,28,40,10,22,34) & (5,17,29,41,11,23,35) &  \\
\hline
(0,3,33,1,21,20,25) & (2,38,18,9,36,26,13) &  & \multirow{2}{*}{$\left<\alpha^{7}\right>$} \\
(0,31,23,19,40,41,32) & (1,6,24,16,25,22,35) &  &  \\
\hline
(0,19,41,17,28,39,4) & (1,5,6,22,37,40,26) & (2,32,15,9,38,36,29) & \multirow{2}{*}{$1$} \\
(3,12,14,35,21,34,18) & (7,30,25,20,24,31,8) & (10,33,23,16,13,11,27) &  \\
\hline
\end{array}}\]
\end{exb}

\begin{exb}\label{hw42r4s37} $\mathrm{HW}^*(42;3^{4},7^{37})$. There are 2 orbits of length 21 and 2 orbits of length 7 on 3-cycles, and 10 orbits of length 21 and 4 orbits of length 3 on 7-cycles.
\setlength{\arraycolsep}{3pt}
\[\fitarray{\begin{array}{|l l l |c|}
\hline
(0,14,28) & (1,29,15) &  & C_{21} \\
\hline
(0,40,15) & (1,8,35) &  & \left<\alpha^{3}\right> \\
\hline
(0,24,6,30,12,36,18) & (1,31,19,7,37,25,13) &  & \multirow{2}{*}{$C_{21}$} \\
(0,18,36,12,30,6,24) & (1,13,25,37,7,19,31) &  &  \\
\hline
(0,34,8,29,10,33,27) & (2,35,39,18,12,9,31) &  & \multirow{2}{*}{$\left<\alpha^{7}\right>$} \\
(0,1,41,31,30,23,12) & (4,11,6,38,21,5,8) &  &  \\
\hline
(0,4,9,29,22,25,34) & (1,35,10,12,32,21,37) & (2,39,11,36,20,19,15) & \multirow{2}{*}{$1$} \\
(3,13,31,33,38,5,16) & (6,17,41,18,28,40,26) & (7,8,23,14,27,24,30) &  \\
\hline
\end{array}}\]
\end{exb}

\begin{exb}\label{hw42r5s36} $\mathrm{HW}^*(42;3^{5},7^{36})$. There are 2 orbits of length 21 and 4 orbits of length 7 on 3-cycles, and 10 orbits of length 21 and 2 orbits of length 3 on 7-cycles.
\setlength{\arraycolsep}{3pt}
\[\fitarray{\begin{array}{|l l l |c|}
\hline
(0,14,28) & (1,29,15) &  & \multirow{2}{*}{$C_{21}$} \\
(0,28,14) & (1,15,29) &  &  \\
\hline
(0,34,17) & (1,38,3) &  & \left<\alpha^{3}\right> \\
\hline
(0,24,6,30,12,36,18) & (1,19,37,13,31,7,25) &  & C_{21} \\
\hline
(0,5,11,41,31,22,16) & (1,4,37,21,20,10,40) &  & \multirow{2}{*}{$\left<\alpha^{7}\right>$} \\
(0,23,5,10,26,11,32) & (1,8,21,41,34,31,2) &  &  \\
\hline
(0,6,10,19,34,11,3) & (1,18,16,33,35,15,31) & (2,4,41,26,37,14,22) & \multirow{2}{*}{$1$} \\
(5,28,24,13,7,36,9) & (8,20,21,30,23,12,32) & (17,25,29,39,40,27,38) &  \\
\hline
\end{array}}\]
\end{exb}

\begin{exb}\label{hw42r6s35} $\mathrm{HW}^*(42;3^{6},7^{35})$. There are 4 orbits of length 21 on 3-cycles, and 10 orbits of length 21 on 7-cycles.
\setlength{\arraycolsep}{3pt}
\[\fitarray{\begin{array}{|l l l |c|}
\hline
(0,32,1) & (3,34,29) &  & \multirow{2}{*}{$\left<\alpha^{3}\right>$} \\
(0,21,23) & (1,26,10) &  &  \\
\hline
(0,20,1,18,23,33,40) & (2,24,41,3,8,21,11) &  & \multirow{2}{*}{$\left<\alpha^{7}\right>$} \\
(0,27,6,37,38,4,8) & (1,35,5,16,25,40,17) &  &  \\
\hline
(0,38,36,24,11,7,29) & (1,41,26,12,6,13,19) & (2,27,20,17,35,21,30) & \multirow{2}{*}{$1$} \\
(3,23,10,16,40,14,15) & (4,39,5,8,18,33,28) & (9,32,31,25,22,34,37) &  \\
\hline
\end{array}}\]
\end{exb}

\begin{exb}\label{hw42r7s34} $\mathrm{HW}^*(42;3^{7},7^{34})$. There are 4 orbits of length 21 and 2 orbits of length 7 on 3-cycles, and 8 orbits of length 21 and 12 orbits of length 3 on 7-cycles.
\setlength{\arraycolsep}{3pt}
\[\fitarray{\begin{array}{|l l l |c|}
\hline
(0,28,14) & (1,15,29) & (2,36,3) & \multirow{2}{*}{$\left<\alpha^{7}\right>$} \\
(4,39,13) & (5,37,20) & (7,24,26) &  \\
\hline
(0,24,6,30,12,36,18) & (1,37,31,25,19,13,7) &  & \multirow{6}{*}{$C_{21}$} \\
(0,6,12,18,24,30,36) & (1,31,19,7,37,25,13) &  &  \\
(0,12,24,36,6,18,30) & (1,19,37,13,31,7,25) &  &  \\
(0,36,30,24,18,12,6) & (1,13,25,37,7,19,31) &  &  \\
(0,18,36,12,30,6,24) & (1,25,7,31,13,37,19) &  &  \\
(0,30,18,6,36,24,12) & (1,7,13,19,25,31,37) &  &  \\
\hline
(0,26,23,7,4,2,15) & (3,8,41,33,20,10,25) &  & \left<\alpha^{7}\right> \\
\hline
(0,41,8,28,17,12,26) & (1,23,21,34,9,13,24) & (2,3,11,30,35,37,40) & \multirow{2}{*}{$1$} \\
(4,7,38,33,29,10,31) & (5,15,22,32,39,25,18) & (6,14,36,19,20,16,27) &  \\
\hline
\end{array}}\]
\end{exb}

\begin{exb}\label{hw42r8s33} $\mathrm{HW}^*(42;3^{8},7^{33})$. There are 4 orbits of length 21 and 4 orbits of length 7 on 3-cycles, and 8 orbits of length 21 and 10 orbits of length 3 on 7-cycles.
\[\fitarray{\begin{array}{|l l |c|}
\hline
(0,28,14) & (1,29,15) & \multirow{2}{*}{$C_{21}$} \\
(0,14,28) & (1,15,29) &  \\
\hline
(0,40,21) & (1,2,17) & \multirow{2}{*}{$\left<\alpha^{3}\right>$} \\
(0,22,23) & (1,26,3) &  \\
\hline
(0,12,24,36,6,18,30) & (1,31,19,7,37,25,13) & \multirow{5}{*}{$C_{21}$} \\
(0,36,30,24,18,12,6) & (1,19,37,13,31,7,25) &  \\
(0,18,36,12,30,6,24) & (1,7,13,19,25,31,37) &  \\
(0,6,12,18,24,30,36) & (1,37,31,25,19,13,7) &  \\
(0,24,6,30,12,36,18) & (1,25,7,31,13,37,19) &  \\
\hline
(0,33,6,13,26,31,1) & (2,36,21,24,23,18,39) & \multirow{4}{*}{$\left<\alpha^{7}\right>$} \\
(0,30,39,19,10,3,35) & (1,9,13,40,6,22,32) &  \\
(0,38,22,9,32,34,31) & (1,39,7,16,41,33,26) &  \\
(0,31,33,11,8,21,10) & (1,30,34,37,12,32,27) &  \\
\hline
\end{array}}\]
\end{exb}

\begin{exb}\label{hw42r9s32} $\mathrm{HW}^*(42;3^{9},7^{32})$. There are 6 orbits of length 21 on 3-cycles, and 8 orbits of length 21 and 8 orbits of length 3 on 7-cycles.
\setlength{\arraycolsep}{3pt}
\[\fitarray{\begin{array}{|l l l |c|}
\hline
(0,35,28) & (1,9,20) &  & \multirow{3}{*}{$\left<\alpha^{3}\right>$} \\
(0,40,25) & (2,33,11) &  &  \\
(0,15,11) & (1,38,22) &  &  \\
\hline
(0,12,24,36,6,18,30) & (1,13,25,37,7,19,31) &  & C_{21} \\
\hline
(0,24,6,30,12,36,18) & (1,7,13,19,25,31,37) & (2,20,38,14,32,8,26) & \multirow{2}{*}{$\left<\alpha^{3}\right>$} \\
(3,33,21,9,39,27,15) & (4,10,16,22,28,34,40) & (5,23,41,17,35,11,29) &  \\
\hline
(0,1,35,2,4,36,23) & (3,6,11,38,5,27,26) &  & \left<\alpha^{7}\right> \\
\hline
(0,30,25,12,34,38,4) & (1,37,41,9,24,21,3) & (2,27,10,13,14,5,31) & \multirow{2}{*}{$1$} \\
(6,26,20,33,19,40,32) & (7,23,28,35,16,15,17) & (8,18,29,36,22,39,11) &  \\
\hline
\end{array}}\]
\end{exb}

\begin{exb}\label{hw42r10s31} $\mathrm{HW}^*(42;3^{10},7^{31})$. There are 6 orbits of length 21 and 2 orbits of length 7 on 3-cycles, and 8 orbits of length 21 and 6 orbits of length 3 on 7-cycles.
\setlength{\arraycolsep}{3pt}
\[\fitarray{\begin{array}{|l l l |c|}
\hline
(0,3,38) & (1,4,23) &  & \left<\alpha^{3}\right> \\
\hline
(0,14,28) & (1,29,15) & (2,11,37) & \multirow{2}{*}{$\left<\alpha^{7}\right>$} \\
(3,27,10) & (4,33,22) & (6,21,40) &  \\
\hline
(0,30,18,6,36,24,12) & (1,37,31,25,19,13,7) & (2,20,38,14,32,8,26) & \multirow{2}{*}{$\left<\alpha^{3}\right>$} \\
(3,21,39,15,33,9,27) & (4,16,28,40,10,22,34) & (5,17,29,41,11,23,35) &  \\
\hline
(0,22,12,37,25,29,3) & (2,18,19,34,41,10,7) &  & \left<\alpha^{7}\right> \\
\hline
(0,17,34,36,5,38,33) & (1,15,21,2,35,27,12) & (3,32,26,25,4,24,37) & \multirow{2}{*}{$1$} \\
(6,40,19,29,9,22,20) & (7,8,18,23,28,13,11) & (10,16,39,41,31,30,14) &  \\
\hline
\end{array}}\]
\end{exb}

\begin{exb}\label{hw42r11s30} $\mathrm{HW}^*(42;3^{11},7^{30})$. There are 6 orbits of length 21 and 4 orbits of length 7 on 3-cycles, and 8 orbits of length 21 and 4 orbits of length 3 on 7-cycles.
\[\fitarray{\begin{array}{|l l |c|}
\hline
(0,28,14) & (1,29,15) & \multirow{2}{*}{$C_{21}$} \\
(0,14,28) & (1,15,29) &  \\
\hline
(0,40,11) & (1,27,26) & \multirow{3}{*}{$\left<\alpha^{3}\right>$} \\
(0,34,21) & (1,8,11) &  \\
(0,26,17) & (1,21,16) &  \\
\hline
(0,24,6,30,12,36,18) & (1,37,31,25,19,13,7) & \multirow{2}{*}{$C_{21}$} \\
(0,18,36,12,30,6,24) & (1,7,13,19,25,31,37) &  \\
\hline
(0,9,18,33,26,7,15) & (2,8,24,13,11,34,31) & \multirow{4}{*}{$\left<\alpha^{7}\right>$} \\
(0,2,10,21,12,33,6) & (1,17,9,27,32,25,36) &  \\
(0,32,15,25,37,19,41) & (2,7,3,22,26,38,6) &  \\
(0,1,5,8,3,32,20) & (2,21,23,40,39,27,24) &  \\
\hline
\end{array}}\]
\end{exb}

\begin{exb}\label{hw42r12s29} $\mathrm{HW}^*(42;3^{12},7^{29})$. There are 8 orbits of length 21 on 3-cycles, and 8 orbits of length 21 and 2 orbits of length 3 on 7-cycles.
\setlength{\arraycolsep}{3pt}
\[\fitarray{\begin{array}{|l l l |c|}
\hline
(0,8,37) & (3,5,28) &  & \multirow{4}{*}{$\left<\alpha^{3}\right>$} \\
(0,38,3) & (1,35,22) &  &  \\
(0,16,5) & (1,39,32) &  &  \\
(0,20,39) & (1,23,34) &  &  \\
\hline
(0,24,6,30,12,36,18) & (1,19,37,13,31,7,25) &  & C_{21} \\
\hline
(0,2,38,19,40,8,35) & (1,25,41,32,31,20,23) &  & \left<\alpha^{7}\right> \\
\hline
(0,14,27,39,10,35,1) & (2,28,21,11,41,16,38) & (3,7,8,9,19,13,40) & \multirow{2}{*}{$1$} \\
(4,22,17,36,26,24,12) & (5,31,37,20,6,18,33) & (15,30,34,25,23,32,29) &  \\
\hline
\end{array}}\]
\end{exb}

\begin{exb}\label{hw42r13s28} $\mathrm{HW}^*(42;3^{13},7^{28})$. There are 8 orbits of length 21 and 2 orbits of length 7 on 3-cycles, and 8 orbits of length 21 on 7-cycles.
\setlength{\arraycolsep}{3pt}
\[\fitarray{\begin{array}{|l l l |c|}
\hline
(0,14,28) & (1,29,15) &  & C_{21} \\
\hline
(0,34,17) & (1,33,8) &  & \multirow{4}{*}{$\left<\alpha^{3}\right>$} \\
(0,40,39) & (1,2,29) &  &  \\
(0,23,20) & (1,9,40) &  &  \\
(0,17,37) & (2,9,16) &  &  \\
\hline
(0,33,18,20,15,39,7) & (2,31,23,36,40,13,24) &  & \left<\alpha^{7}\right> \\
\hline
(0,19,40,8,39,35,12) & (1,10,36,33,28,22,34) & (2,15,9,21,27,18,29) & \multirow{2}{*}{$1$} \\
(3,26,32,37,41,25,5) & (4,20,16,24,14,38,17) & (6,7,23,11,13,31,30) &  \\
\hline
\end{array}}\]
\end{exb}

\begin{exb}\label{hw42r14s27} $\mathrm{HW}^*(42;3^{14},7^{27})$. There are 8 orbits of length 21 and 4 orbits of length 7 on 3-cycles, and 6 orbits of length 21 and 12 orbits of length 3 on 7-cycles.
\setlength{\arraycolsep}{3pt}
\[\fitarray{\begin{array}{|l l l |c|}
\hline
(0,28,14) & (1,29,15) & (2,12,37) & \multirow{2}{*}{$\left<\alpha^{7}\right>$} \\
(3,25,8) & (4,27,38) & (5,20,7) &  \\
\hline
(0,14,28) & (1,15,29) & (2,34,13) & \multirow{2}{*}{$\left<\alpha^{7}\right>$} \\
(3,12,4) & (5,37,38) & (7,36,11) &  \\
\hline
(0,12,24,36,6,18,30) & (1,19,37,13,31,7,25) &  & \multirow{3}{*}{$C_{21}$} \\
(0,6,12,18,24,30,36) & (1,31,19,7,37,25,13) &  &  \\
(0,30,18,6,36,24,12) & (1,25,7,31,13,37,19) &  &  \\
\hline
(0,18,36,12,30,6,24) & (1,7,13,19,25,31,37) & (2,38,32,26,20,14,8) & \multirow{2}{*}{$\left<\alpha^{3}\right>$} \\
(3,15,27,39,9,21,33) & (4,28,10,34,16,40,22) & (5,41,35,29,23,17,11) &  \\
\hline
(0,33,18,31,22,27,20) & (1,9,11,14,21,38,12) & (2,3,16,35,39,13,5) & \multirow{2}{*}{$1$} \\
(4,7,6,10,30,41,25) & (8,23,28,26,19,29,24) & (15,34,36,32,17,40,37) &  \\
\hline
\end{array}}\]
\end{exb}

\begin{exb}\label{hw42r15s26} $\mathrm{HW}^*(42;3^{15},7^{26})$. There are 10 orbits of length 21 on 3-cycles, and 6 orbits of length 21 and 10 orbits of length 3 on 7-cycles.
\setlength{\arraycolsep}{3pt}
\[\fitarray{\begin{array}{|l l l |c|}
\hline
(0,33,13) & (2,5,22) &  & \multirow{5}{*}{$\left<\alpha^{3}\right>$} \\
(0,40,19) & (2,33,5) &  &  \\
(0,26,22) & (1,33,35) &  &  \\
(0,34,9) & (1,11,32) &  &  \\
(0,10,23) & (1,39,38) &  &  \\
\hline
(0,30,18,6,36,24,12) & (1,37,31,25,19,13,7) &  & \multirow{5}{*}{$C_{21}$} \\
(0,24,6,30,12,36,18) & (1,7,13,19,25,31,37) &  &  \\
(0,6,12,18,24,30,36) & (1,25,7,31,13,37,19) &  &  \\
(0,36,30,24,18,12,6) & (1,19,37,13,31,7,25) &  &  \\
(0,18,36,12,30,6,24) & (1,13,25,37,7,19,31) &  &  \\
\hline
(0,29,30,7,32,27,11) & (1,8,15,19,35,21,9) & (2,34,6,22,26,41,39) & \multirow{2}{*}{$1$} \\
(3,23,38,31,16,28,36) & (4,5,14,13,24,10,33) & (12,37,40,25,18,20,17) &  \\
\hline
\end{array}}\]
\end{exb}

\begin{exb}\label{hw42r16s25} $\mathrm{HW}^*(42;3^{16},7^{25})$. There are 10 orbits of length 21 and 2 orbits of length 7 on 3-cycles, and 6 orbits of length 21 and 8 orbits of length 3 on 7-cycles.
\setlength{\arraycolsep}{3pt}
\[\fitarray{\begin{array}{|l l l |c|}
\hline
(0,14,28) & (1,29,15) &  & C_{21} \\
\hline
(0,40,21) & (1,32,17) &  & \multirow{5}{*}{$\left<\alpha^{3}\right>$} \\
(0,21,37) & (2,40,5) &  &  \\
(0,31,9) & (2,34,35) &  &  \\
(0,5,27) & (1,26,4) &  &  \\
(0,16,20) & (1,9,41) &  &  \\
\hline
(0,6,12,18,24,30,36) & (1,13,25,37,7,19,31) &  & C_{21} \\
\hline
(0,12,24,36,6,18,30) & (1,25,7,31,13,37,19) & (2,32,20,8,38,26,14) & \multirow{2}{*}{$\left<\alpha^{3}\right>$} \\
(3,21,39,15,33,9,27) & (4,22,40,16,34,10,28) & (5,11,17,23,29,35,41) &  \\
\hline
(0,26,25,13,16,3,1) & (2,30,22,4,17,9,36) & (5,19,12,29,10,21,32) & \multirow{2}{*}{$1$} \\
(6,41,35,39,34,37,38) & (7,24,18,20,15,28,11) & (8,27,14,23,33,40,31) &  \\
\hline
\end{array}}\]
\end{exb}

\begin{exb}\label{hw42r17s24} $\mathrm{HW}^*(42;3^{17},7^{24})$. There are 10 orbits of length 21 and 4 orbits of length 7 on 3-cycles, and 6 orbits of length 21 and 6 orbits of length 3 on 7-cycles.
\setlength{\arraycolsep}{3pt}
\[\fitarray{\begin{array}{|l l l |c|}
\hline
(0,5,7) & (2,21,34) &  & \left<\alpha^{3}\right> \\
\hline
(0,14,28) & (1,15,29) & (2,41,22) & \multirow{2}{*}{$\left<\alpha^{7}\right>$} \\
(3,24,19) & (4,39,26) & (6,21,9) &  \\
\hline
(0,28,14) & (1,29,15) & (2,33,6) & \multirow{2}{*}{$\left<\alpha^{7}\right>$} \\
(3,22,35) & (4,27,12) & (9,10,11) &  \\
\hline
(0,6,12,18,24,30,36) & (1,37,31,25,19,13,7) & (2,38,32,26,20,14,8) & \multirow{2}{*}{$\left<\alpha^{3}\right>$} \\
(3,27,9,33,15,39,21) & (4,22,40,16,34,10,28) & (5,11,17,23,29,35,41) &  \\
\hline
(0,40,30,41,32,20,5) & (1,8,10,39,9,31,35) &  & \multirow{3}{*}{$\left<\alpha^{7}\right>$} \\
(0,41,10,22,5,21,17) & (1,4,25,30,34,37,26) &  &  \\
(0,26,8,17,16,7,25) & (1,10,18,34,41,19,9) &  &  \\
\hline
\end{array}}\]
\end{exb}

\begin{exb}\label{hw42r18s23} $\mathrm{HW}^*(42;3^{18},7^{23})$. There are 12 orbits of length 21 on 3-cycles, and 6 orbits of length 21 and 4 orbits of length 3 on 7-cycles.
\setlength{\arraycolsep}{3pt}
\[\fitarray{\begin{array}{|l l l |c|}
\hline
(0,11,31) & (2,9,10) &  & \multirow{6}{*}{$\left<\alpha^{3}\right>$} \\
(0,32,7) & (3,22,23) &  &  \\
(0,35,21) & (1,26,10) &  &  \\
(0,2,15) & (1,34,11) &  &  \\
(0,3,29) & (1,10,38) &  &  \\
(0,25,28) & (2,41,15) &  &  \\
\hline
(0,24,6,30,12,36,18) & (1,37,31,25,19,13,7) &  & \multirow{2}{*}{$C_{21}$} \\
(0,12,24,36,6,18,30) & (1,13,25,37,7,19,31) &  &  \\
\hline
(0,37,29,24,22,16,4) & (1,31,12,41,40,25,39) & (2,23,5,15,21,19,27) & \multirow{2}{*}{$1$} \\
(3,34,33,35,11,18,7) & (6,28,36,9,14,20,38) & (8,17,32,13,10,26,30) &  \\
\hline
\end{array}}\]
\end{exb}

\begin{exb}\label{hw42r19s22} $\mathrm{HW}^*(42;3^{19},7^{22})$. There are 12 orbits of length 21 and 2 orbits of length 7 on 3-cycles, and 6 orbits of length 21 and 2 orbits of length 3 on 7-cycles.
\[\fitarray{\begin{array}{|l l |c|}
\hline
(0,14,28) & (1,29,15) & C_{21} \\
\hline
(0,32,3) & (1,23,10) & \multirow{6}{*}{$\left<\alpha^{3}\right>$} \\
(0,40,8) & (1,27,5) &  \\
(0,38,37) & (3,17,34) &  \\
(0,28,29) & (1,20,27) &  \\
(0,3,26) & (1,35,28) &  \\
(0,37,40) & (2,33,23) &  \\
\hline
(0,12,24,36,6,18,30) & (1,31,19,7,37,25,13) & C_{21} \\
\hline
(0,23,31,4,10,7,16) & (1,25,20,40,22,27,33) & \multirow{3}{*}{$\left<\alpha^{7}\right>$} \\
(0,17,35,5,16,34,1) & (4,8,38,23,25,26,13) &  \\
(0,8,27,34,17,11,15) & (2,37,35,24,4,40,19) &  \\
\hline
\end{array}}\]
\end{exb}

\begin{exb}\label{hw42r20s21} $\mathrm{HW}^*(42;3^{20},7^{21})$. There are 12 orbits of length 21 and 4 orbits of length 7 on 3-cycles, and 6 orbits of length 21 on 7-cycles.
\[\fitarray{\begin{array}{|l l l |c|}
\hline
(0,34,31) & (2,35,15) &  & \multirow{2}{*}{$\left<\alpha^{3}\right>$} \\
(0,37,5) & (2,10,39) &  &  \\
\hline
(0,14,28) & (1,15,29) & (2,24,5) & \multirow{2}{*}{$\left<\alpha^{7}\right>$} \\
(3,39,27) & (4,20,26) & (7,23,8) &  \\
\hline
(0,28,14) & (1,29,15) & (2,4,8) & \multirow{2}{*}{$\left<\alpha^{7}\right>$} \\
(3,34,5) & (7,26,41) & (9,11,24) &  \\
\hline
(0,3,10,41,20,2,19) & (1,39,21,12,37,4,22) &  & \multirow{3}{*}{$\left<\alpha^{7}\right>$} \\
(0,1,21,33,37,27,2) & (3,18,11,10,6,36,26) &  &  \\
(0,11,37,38,8,15,7) & (2,12,17,20,4,13,19) &  &  \\
\hline
\end{array}}\]
\end{exb}

\begin{exb}\label{hw42r21s20} $\mathrm{HW}^*(42;3^{21},7^{20})$. There are 14 orbits of length 21 on 3-cycles, and 4 orbits of length 21 and 12 orbits of length 3 on 7-cycles.
\setlength{\arraycolsep}{3pt}
\[\fitarray{\begin{array}{|l l l |c|}
\hline
(0,28,33) & (1,15,14) & (2,37,41) & \multirow{5}{*}{$1$} \\
(3,32,17) & (4,23,12) & (5,26,24) &  \\
(6,8,22) & (7,20,16) & (9,29,13) &  \\
(10,35,27) & (11,34,38) & (18,31,21) &  \\
(19,36,39) & (25,40,30) &  &  \\
\hline
(0,12,24,36,6,18,30) & (1,31,19,7,37,25,13) & (2,38,32,26,20,14,8) & \multirow{2}{*}{$\left<\alpha^{3}\right>$} \\
(3,15,27,39,9,21,33) & (4,10,16,22,28,34,40) & (5,11,17,23,29,35,41) &  \\
\hline
(0,24,6,30,12,36,18) & (1,37,31,25,19,13,7) & (2,32,20,8,38,26,14) & \multirow{2}{*}{$\left<\alpha^{3}\right>$} \\
(3,27,9,33,15,39,21) & (4,22,40,16,34,10,28) & (5,23,41,17,35,11,29) &  \\
\hline
(0,21,40,5,13,20,9) & (1,3,8,30,39,24,32) &  & \multirow{2}{*}{$\left<\alpha^{7}\right>$} \\
(0,10,27,25,20,21,22) & (1,12,9,2,18,17,33) &  &  \\
\hline
\end{array}}\]
\end{exb}

\begin{exb}\label{hw42r22s19} $\mathrm{HW}^*(42;3^{22},7^{19})$. There are 14 orbits of length 21 and 2 orbits of length 7 on 3-cycles, and 4 orbits of length 21 and 10 orbits of length 3 on 7-cycles.
\[\fitarray{\begin{array}{|l l |c|}
\hline
(0,28,14) & (1,15,29) & C_{21} \\
\hline
(0,10,8) & (1,29,9) & \multirow{7}{*}{$\left<\alpha^{3}\right>$} \\
(0,8,7) & (3,29,28) &  \\
(0,14,23) & (1,21,16) &  \\
(0,20,27) & (1,17,40) &  \\
(0,19,17) & (2,39,28) &  \\
(0,31,21) & (2,35,22) &  \\
(0,39,40) & (1,5,32) &  \\
\hline
(0,30,18,6,36,24,12) & (1,25,7,31,13,37,19) & \multirow{5}{*}{$C_{21}$} \\
(0,36,30,24,18,12,6) & (1,13,25,37,7,19,31) &  \\
(0,18,36,12,30,6,24) & (1,37,31,25,19,13,7) &  \\
(0,6,12,18,24,30,36) & (1,31,19,7,37,25,13) &  \\
(0,12,24,36,6,18,30) & (1,7,13,19,25,31,37) &  \\
\hline
(0,23,33,24,25,34,18) & (1,12,27,30,17,22,35) & \multirow{2}{*}{$\left<\alpha^{7}\right>$} \\
(0,25,1,8,12,33,4) & (2,7,9,6,38,31,27) &  \\
\hline
\end{array}}\]
\end{exb}

\begin{exb}\label{hw42r23s18} $\mathrm{HW}^*(42;3^{23},7^{18})$. There are 14 orbits of length 21 and 4 orbits of length 7 on 3-cycles, and 4 orbits of length 21 and 8 orbits of length 3 on 7-cycles.
\setlength{\arraycolsep}{3pt}
\[\fitarray{\begin{array}{|l l l |c|}
\hline
(0,14,28) & (1,15,29) &  & \multirow{2}{*}{$C_{21}$} \\
(0,28,14) & (1,29,15) &  &  \\
\hline
(0,10,2) & (1,37,6) & (3,32,29) & \multirow{5}{*}{$1$} \\
(4,17,21) & (5,36,12) & (7,34,38) &  \\
(8,27,35) & (9,11,18) & (13,26,33) &  \\
(14,19,39) & (15,24,23) & (16,41,20) &  \\
(22,30,31) & (25,28,40) &  &  \\
\hline
(0,6,12,18,24,30,36) & (1,31,19,7,37,25,13) &  & C_{21} \\
\hline
(0,30,18,6,36,24,12) & (1,25,7,31,13,37,19) & (2,26,8,32,14,38,20) & \multirow{2}{*}{$\left<\alpha^{3}\right>$} \\
(3,9,15,21,27,33,39) & (4,40,34,28,22,16,10) & (5,17,29,41,11,23,35) &  \\
\hline
(0,17,1,19,24,2,23) & (4,27,22,11,12,34,7) &  & \multirow{2}{*}{$\left<\alpha^{7}\right>$} \\
(0,32,19,15,25,23,16) & (3,26,35,34,36,10,13) &  &  \\
\hline
\end{array}}\]
\end{exb}

\begin{exb}\label{hw42r24s17} $\mathrm{HW}^*(42;3^{24},7^{17})$. There are 16 orbits of length 21 on 3-cycles, and 4 orbits of length 21 and 6 orbits of length 3 on 7-cycles.
\setlength{\arraycolsep}{3pt}
\[\fitarray{\begin{array}{|l l l |c|}
\hline
(0,17,13) & (2,21,28) &  & \multirow{8}{*}{$\left<\alpha^{3}\right>$} \\
(0,1,4) & (2,39,17) &  &  \\
(0,39,29) & (1,2,16) &  &  \\
(0,31,8) & (3,5,22) &  &  \\
(0,32,40) & (1,23,15) &  &  \\
(0,29,33) & (1,32,10) &  &  \\
(0,13,21) & (2,17,16) &  &  \\
(0,22,27) & (1,26,29) &  &  \\
\hline
(0,18,36,12,30,6,24) & (1,37,31,25,19,13,7) & (2,38,32,26,20,14,8) & \multirow{2}{*}{$\left<\alpha^{3}\right>$} \\
(3,15,27,39,9,21,33) & (4,10,16,22,28,34,40) & (5,23,41,17,35,11,29) &  \\
\hline
(0,24,8,7,37,5,31) & (1,40,2,27,32,39,34) &  & \multirow{2}{*}{$\left<\alpha^{7}\right>$} \\
(0,10,22,33,15,31,12) & (2,37,35,41,34,32,11) &  &  \\
\hline
\end{array}}\]
\end{exb}

\begin{exb}\label{hw42r25s16} $\mathrm{HW}^*(42;3^{25},7^{16})$. There are 16 orbits of length 21 and 2 orbits of length 7 on 3-cycles, and 4 orbits of length 21 and 4 orbits of length 3 on 7-cycles.
\setlength{\arraycolsep}{3pt}
\[\fitarray{\begin{array}{|l l l |c|}
\hline
(0,28,37) & (2,35,3) &  & \multirow{6}{*}{$\left<\alpha^{3}\right>$} \\
(0,38,17) & (1,22,9) &  &  \\
(0,31,39) & (2,22,5) &  &  \\
(0,27,34) & (1,20,17) &  &  \\
(0,11,40) & (1,15,38) &  &  \\
(0,19,8) & (3,35,28) &  &  \\
\hline
(0,14,28) & (1,29,15) & (2,20,35) & \multirow{2}{*}{$\left<\alpha^{7}\right>$} \\
(3,5,41) & (4,39,9) & (8,12,24) &  \\
\hline
(0,6,12,18,24,30,36) & (1,31,19,7,37,25,13) &  & \multirow{2}{*}{$C_{21}$} \\
(0,30,18,6,36,24,12) & (1,19,37,13,31,7,25) &  &  \\
\hline
(0,7,18,34,16,26,27) & (1,25,10,9,5,22,3) &  & \multirow{2}{*}{$\left<\alpha^{7}\right>$} \\
(0,13,33,4,7,29,20) & (2,38,36,31,37,11,12) &  &  \\
\hline
\end{array}}\]
\end{exb}

\begin{exb}\label{hw42r26s15} $\mathrm{HW}^*(42;3^{26},7^{15})$. There are 16 orbits of length 21 and 4 orbits of length 7 on 3-cycles, and 4 orbits of length 21 and 2 orbits of length 3 on 7-cycles.
\setlength{\arraycolsep}{3pt}
\[\fitarray{\begin{array}{|l l l |c|}
\hline
(0,14,28) & (1,29,15) &  & \multirow{2}{*}{$C_{21}$} \\
(0,28,14) & (1,15,29) &  &  \\
\hline
(0,16,41) & (1,39,20) &  & \left<\alpha^{3}\right> \\
\hline
(0,11,6) & (1,4,13) & (2,17,36) & \multirow{5}{*}{$1$} \\
(3,10,37) & (5,29,14) & (7,18,15) &  \\
(8,30,32) & (9,34,35) & (12,19,23) &  \\
(16,33,26) & (20,39,41) & (21,38,31) &  \\
(22,25,40) & (24,28,27) &  &  \\
\hline
(0,6,12,18,24,30,36) & (1,37,31,25,19,13,7) &  & C_{21} \\
\hline
(0,26,38,33,13,31,22) & (1,6,4,25,9,21,30) &  & \multirow{2}{*}{$\left<\alpha^{7}\right>$} \\
(0,10,6,40,27,33,1) & (2,32,37,35,22,11,31) &  &  \\
\hline
\end{array}}\]
\end{exb}

\begin{exb}\label{hw42r27s14} $\mathrm{HW}^*(42;3^{27},7^{14})$. There are 18 orbits of length 21 on 3-cycles, and 4 orbits of length 21 on 7-cycles.
\setlength{\arraycolsep}{3pt}
\[\fitarray{\begin{array}{|l l l |c|}
\hline
(0,14,1) & (3,29,10) &  & \multirow{2}{*}{$\left<\alpha^{3}\right>$} \\
(0,15,22) & (1,23,2) &  &  \\
\hline
(0,22,17) & (1,33,19) & (2,21,15) & \multirow{5}{*}{$1$} \\
(3,14,35) & (4,27,18) & (5,23,39) &  \\
(6,40,9) & (7,8,10) & (11,41,13) &  \\
(12,20,36) & (16,25,28) & (24,29,31) &  \\
(26,30,37) & (32,38,34) &  &  \\
\hline
(0,13,40,24,7,3,37) & (1,5,25,20,32,22,16) &  & \multirow{2}{*}{$\left<\alpha^{7}\right>$} \\
(0,17,29,38,20,21,27) & (2,5,18,9,26,36,25) &  &  \\
\hline
\end{array}}\]
\end{exb}

\begin{exb}\label{hw42r28s13} $\mathrm{HW}^*(42;3^{28},7^{13})$. There are 18 orbits of length 21 and 2 orbits of length 7 on 3-cycles, and 2 orbits of length 21 and 12 orbits of length 3 on 7-cycles.
\setlength{\arraycolsep}{3pt}
\[\fitarray{\begin{array}{|l l l |c|}
\hline
(0,14,28) & (1,29,15) & (2,34,21) & \multirow{2}{*}{$\left<\alpha^{7}\right>$} \\
(3,38,11) & (4,9,36) & (5,13,26) &  \\
\hline
(0,35,9) & (1,2,5) & (3,12,34) & \multirow{5}{*}{$1$} \\
(4,17,15) & (6,37,20) & (7,29,26) &  \\
(8,28,24) & (10,19,33) & (11,31,36) &  \\
(13,39,18) & (14,16,41) & (21,23,22) &  \\
(25,32,40) & (27,38,30) &  &  \\
\hline
(0,30,18,6,36,24,12) & (1,25,7,31,13,37,19) & (2,8,14,20,26,32,38) & \multirow{2}{*}{$\left<\alpha^{3}\right>$} \\
(3,33,21,9,39,27,15) & (4,22,40,16,34,10,28) & (5,23,41,17,35,11,29) &  \\
\hline
(0,12,24,36,6,18,30) & (1,13,25,37,7,19,31) & (2,26,8,32,14,38,20) & \multirow{2}{*}{$\left<\alpha^{3}\right>$} \\
(3,9,15,21,27,33,39) & (4,40,34,28,22,16,10) & (5,41,35,29,23,17,11) &  \\
\hline
(0,19,29,32,36,34,25) & (2,9,38,12,13,3,7) &  & \left<\alpha^{7}\right> \\
\hline
\end{array}}\]
\end{exb}

\begin{exb}\label{hw42r29s12} $\mathrm{HW}^*(42;3^{29},7^{12})$. There are 18 orbits of length 21 and 4 orbits of length 7 on 3-cycles, and 2 orbits of length 21 and 10 orbits of length 3 on 7-cycles.
\[\fitarray{\begin{array}{|l l l |c|}
\hline
(0,28,14) & (1,15,29) &  & \multirow{2}{*}{$C_{21}$} \\
(0,14,28) & (1,29,15) &  &  \\
\hline
(0,41,26) & (1,22,9) &  & \multirow{2}{*}{$\left<\alpha^{3}\right>$} \\
(0,27,2) & (1,10,5) &  &  \\
\hline
(0,38,19) & (1,16,8) & (2,41,34) & \multirow{5}{*}{$1$} \\
(3,29,24) & (4,9,20) & (5,30,37) &  \\
(6,31,28) & (7,27,40) & (10,18,33) &  \\
(11,14,36) & (12,23,21) & (13,15,39) &  \\
(17,22,26) & (25,32,35) &  &  \\
\hline
(0,12,24,36,6,18,30) & (1,19,37,13,31,7,25) &  & \multirow{5}{*}{$C_{21}$} \\
(0,36,30,24,18,12,6) & (1,7,13,19,25,31,37) &  &  \\
(0,24,6,30,12,36,18) & (1,37,31,25,19,13,7) &  &  \\
(0,30,18,6,36,24,12) & (1,13,25,37,7,19,31) &  &  \\
(0,6,12,18,24,30,36) & (1,31,19,7,37,25,13) &  &  \\
\hline
(0,2,34,5,4,22,41) & (1,9,31,35,24,25,12) &  & \left<\alpha^{7}\right> \\
\hline
\end{array}}\]
\end{exb}

\begin{exb}\label{hw42r30s11} $\mathrm{HW}^*(42;3^{30},7^{11})$. There are 20 orbits of length 21 on 3-cycles, and 2 orbits of length 21 and 8 orbits of length 3 on 7-cycles.
\[\fitarray{\begin{array}{|l l |c|}
\hline
(0,14,27) & (1,40,23) & \multirow{10}{*}{$\left<\alpha^{3}\right>$} \\
(0,8,1) & (3,35,22) &  \\
(0,28,29) & (1,38,33) &  \\
(0,4,20) & (1,3,29) &  \\
(0,40,25) & (2,11,39) &  \\
(0,11,22) & (1,20,27) &  \\
(0,2,35) & (1,10,9) &  \\
(0,19,10) & (2,41,21) &  \\
(0,10,15) & (1,5,26) &  \\
(0,15,11) & (1,26,22) &  \\
\hline
(0,36,30,24,18,12,6) & (1,25,7,31,13,37,19) & \multirow{4}{*}{$C_{21}$} \\
(0,12,24,36,6,18,30) & (1,31,19,7,37,25,13) &  \\
(0,30,18,6,36,24,12) & (1,13,25,37,7,19,31) &  \\
(0,6,12,18,24,30,36) & (1,19,37,13,31,7,25) &  \\
\hline
(0,29,37,38,22,4,7) & (2,33,27,25,31,34,26) & \left<\alpha^{7}\right> \\
\hline
\end{array}}\]
\end{exb}

\begin{exb}\label{hw42r31s10} $\mathrm{HW}^*(42;3^{31},7^{10})$. There are 20 orbits of length 21 and 2 orbits of length 7 on 3-cycles, and 2 orbits of length 21 and 6 orbits of length 3 on 7-cycles.
\setlength{\arraycolsep}{3pt}
\[\fitarray{\begin{array}{|l l l |c|}
\hline
(0,14,28) & (1,29,15) &  & C_{21} \\
\hline
(0,27,29) & (1,2,22) &  & \multirow{3}{*}{$\left<\alpha^{3}\right>$} \\
(0,26,35) & (1,22,33) &  &  \\
(0,41,7) & (2,34,15) &  &  \\
\hline
(0,17,8) & (1,4,37) & (2,5,27) & \multirow{5}{*}{$1$} \\
(3,29,40) & (6,13,33) & (7,34,36) &  \\
(9,14,12) & (10,35,24) & (11,25,15) &  \\
(16,38,20) & (18,28,21) & (19,23,32) &  \\
(22,41,39) & (26,30,31) &  &  \\
\hline
(0,30,18,6,36,24,12) & (1,37,31,25,19,13,7) & (2,20,38,14,32,8,26) & \multirow{2}{*}{$\left<\alpha^{3}\right>$} \\
(3,21,39,15,33,9,27) & (4,40,34,28,22,16,10) & (5,35,23,11,41,29,17) &  \\
\hline
(0,12,27,39,20,9,1) & (2,18,24,19,35,17,36) &  & \left<\alpha^{7}\right> \\
\hline
\end{array}}\]
\end{exb}

\begin{exb}\label{hw42r32s9} $\mathrm{HW}^*(42;3^{32},7^{9})$. There are 20 orbits of length 21 and 4 orbits of length 7 on 3-cycles, and 2 orbits of length 21 and 4 orbits of length 3 on 7-cycles.
\[\fitarray{\begin{array}{|l l l |c|}
\hline
(0,14,28) & (1,29,15) &  & C_{21} \\
\hline
(0,37,34) & (2,21,29) &  & \left<\alpha^{3}\right> \\
\hline
(0,28,14) & (1,15,29) & (2,35,38) & \multirow{2}{*}{$\left<\alpha^{7}\right>$} \\
(3,12,33) & (4,13,11) & (6,23,8) &  \\
\hline
(0,4,11) & (1,12,13) & (2,26,7) & \multirow{5}{*}{$1$} \\
(3,29,25) & (5,15,6) & (8,24,14) &  \\
(9,41,22) & (10,40,32) & (16,19,20) &  \\
(17,30,35) & (18,28,39) & (21,37,36) &  \\
(23,27,34) & (31,33,38) &  &  \\
\hline
(0,18,36,12,30,6,24) & (1,19,37,13,31,7,25) &  & \multirow{2}{*}{$C_{21}$} \\
(0,12,24,36,6,18,30) & (1,37,31,25,19,13,7) &  &  \\
\hline
(0,39,22,2,15,7,13) & (3,20,4,19,38,40,23) &  & \left<\alpha^{7}\right> \\
\hline
\end{array}}\]
\end{exb}

\begin{exb}\label{hw42r33s8} $\mathrm{HW}^*(42;3^{33},7^{8})$. There are 22 orbits of length 21 on 3-cycles, and 2 orbits of length 21 and 2 orbits of length 3 on 7-cycles.
\[\fitarray{\begin{array}{|l l l |c|}
\hline
(0,37,4) & (2,35,15) &  & \multirow{4}{*}{$\left<\alpha^{3}\right>$} \\
(0,17,38) & (1,3,28) &  &  \\
(0,10,37) & (2,27,5) &  &  \\
(0,8,29) & (1,4,9) &  &  \\
\hline
(0,28,9) & (1,19,34) & (2,16,18) & \multirow{5}{*}{$1$} \\
(3,17,10) & (4,36,12) & (5,24,37) &  \\
(6,25,32) & (7,23,11) & (8,15,39) &  \\
(13,30,27) & (14,26,29) & (20,40,38) &  \\
(21,33,22) & (31,35,41) &  &  \\
\hline
(0,6,12,18,24,30,36) & (1,37,31,25,19,13,7) &  & C_{21} \\
\hline
(0,1,27,26,20,8,19) & (2,31,39,37,38,18,7) &  & \left<\alpha^{7}\right> \\
\hline
\end{array}}\]
\end{exb}

\begin{exb}\label{hw42r34s7} $\mathrm{HW}^*(42;3^{34},7^{7})$. There are 22 orbits of length 21 and 2 orbits of length 7 on 3-cycles, and 2 orbits of length 21 on 7-cycles.
\setlength{\arraycolsep}{3pt}
\[\fitarray{\begin{array}{|l l l |c|}
\hline
(0,14,28) & (1,29,15) &  & C_{21} \\
\hline
(0,9,16) & (1,23,8) &  & \multirow{4}{*}{$\left<\alpha^{3}\right>$} \\
(0,4,25) & (2,15,35) &  &  \\
(0,37,34) & (2,41,33) &  &  \\
(0,29,32) & (1,22,27) &  &  \\
\hline
(0,27,18) & (1,9,11) & (2,32,22) & \multirow{5}{*}{$1$} \\
(3,15,34) & (4,40,17) & (5,31,30) &  \\
(6,8,20) & (7,25,26) & (10,16,36) &  \\
(12,37,41) & (13,19,14) & (21,35,33) &  \\
(23,38,29) & (24,39,28) &  &  \\
\hline
(0,7,3,26,27,32,24) & (1,11,36,34,23,16,19) &  & \left<\alpha^{7}\right> \\
\hline
\end{array}}\]
\end{exb}

\begin{exb}\label{hw42r35s6} $\mathrm{HW}^*(42;3^{35},7^{6})$. There are 22 orbits of length 21 and 4 orbits of length 7 on 3-cycles, and 12 orbits of length 3 on 7-cycles.
\setlength{\arraycolsep}{3pt}
\[\fitarray{\begin{array}{|l l l |c|}
\hline
(0,28,14) & (1,15,29) & (2,18,11) & \multirow{2}{*}{$\left<\alpha^{7}\right>$} \\
(3,40,37) & (5,20,22) & (7,24,27) &  \\
\hline
(0,14,28) & (1,29,15) & (2,7,6) & \multirow{2}{*}{$\left<\alpha^{7}\right>$} \\
(3,12,25) & (4,24,41) & (5,9,36) &  \\
\hline
(0,27,29) & (1,24,5) & (2,17,33) & \multirow{5}{*}{$1$} \\
(3,32,16) & (4,23,7) & (6,13,38) &  \\
(8,40,19) & (9,12,20) & (10,11,18) &  \\
(14,36,35) & (15,34,25) & (21,22,26) &  \\
(28,37,30) & (31,41,39) &  &  \\
\hline
(0,36,30,24,18,12,6) & (1,19,37,13,31,7,25) & (2,8,14,20,26,32,38) & \multirow{2}{*}{$\left<\alpha^{3}\right>$} \\
(3,9,15,21,27,33,39) & (4,34,22,10,40,28,16) & (5,41,35,29,23,17,11) &  \\
\hline
(0,24,6,30,12,36,18) & (1,31,19,7,37,25,13) & (2,20,38,14,32,8,26) & \multirow{2}{*}{$\left<\alpha^{3}\right>$} \\
(3,27,9,33,15,39,21) & (4,16,28,40,10,22,34) & (5,17,29,41,11,23,35) &  \\
\hline
\end{array}}\]
\end{exb}

\begin{exb}\label{hw42r36s5} $\mathrm{HW}^*(42;3^{36},7^{5})$. There are 24 orbits of length 21 on 3-cycles, and 10 orbits of length 3 on 7-cycles.
\[\fitarray{\begin{array}{|l l l |c|}
\hline
(0,21,26) & (1,34,23) &  & \multirow{5}{*}{$\left<\alpha^{3}\right>$} \\
(0,28,8) & (1,33,5) &  &  \\
(0,38,1) & (3,28,41) &  &  \\
(0,20,13) & (3,34,11) &  &  \\
(0,14,5) & (1,4,33) &  &  \\
\hline
(0,10,25) & (1,27,34) & (2,13,36) & \multirow{5}{*}{$1$} \\
(3,12,38) & (4,8,31) & (5,7,20) &  \\
(6,18,21) & (9,17,14) & (11,30,28) &  \\
(15,37,16) & (19,35,33) & (22,39,32) &  \\
(23,24,26) & (29,40,41) &  &  \\
\hline
(0,6,12,18,24,30,36) & (1,13,25,37,7,19,31) &  & \multirow{5}{*}{$C_{21}$} \\
(0,24,6,30,12,36,18) & (1,19,37,13,31,7,25) &  &  \\
(0,36,30,24,18,12,6) & (1,25,7,31,13,37,19) &  &  \\
(0,18,36,12,30,6,24) & (1,37,31,25,19,13,7) &  &  \\
(0,30,18,6,36,24,12) & (1,7,13,19,25,31,37) &  &  \\
\hline
\end{array}}\]
\end{exb}

\begin{exb}\label{hw42r37s4} $\mathrm{HW}^*(42;3^{37},7^{4})$. There are 24 orbits of length 21 and 2 orbits of length 7 on 3-cycles, and 8 orbits of length 3 on 7-cycles.
\[\fitarray{\begin{array}{|l l l |c|}
\hline
(0,31,41) & (2,28,3) &  & \multirow{10}{*}{$\left<\alpha^{3}\right>$} \\
(0,27,29) & (1,10,20) &  &  \\
(0,1,5) & (2,22,33) &  &  \\
(0,33,19) & (2,41,28) &  &  \\
(0,25,17) & (2,39,4) &  &  \\
(0,38,15) & (1,22,29) &  &  \\
(0,22,25) & (2,15,5) &  &  \\
(0,32,40) & (1,21,17) &  &  \\
(0,5,38) & (1,16,3) &  &  \\
(0,15,8) & (1,17,22) &  &  \\
\hline
(0,41,30) & (1,20,8) & (3,11,23) & \multirow{2}{*}{$\left<\alpha^{7}\right>$} \\
(4,32,18) & (5,35,38) & (12,26,40) &  \\
\hline
(0,24,6,30,12,36,18) & (1,19,37,13,31,7,25) &  & \multirow{4}{*}{$C_{21}$} \\
(0,6,12,18,24,30,36) & (1,37,31,25,19,13,7) &  &  \\
(0,36,30,24,18,12,6) & (1,25,7,31,13,37,19) &  &  \\
(0,18,36,12,30,6,24) & (1,7,13,19,25,31,37) &  &  \\
\hline
\end{array}}\]
\end{exb}

\begin{exb}\label{hw42r38s3} $\mathrm{HW}^*(42;3^{38},7^{3})$. There are 24 orbits of length 21 and 4 orbits of length 7 on 3-cycles, and 6 orbits of length 3 on 7-cycles.
\setlength{\arraycolsep}{3pt}
\[\fitarray{\begin{array}{|l l l |c|}
\hline
(0,26,33) & (1,17,4) &  & \left<\alpha^{3}\right> \\
\hline
(0,28,14) & (1,15,29) & (2,6,21) & \multirow{2}{*}{$\left<\alpha^{7}\right>$} \\
(3,40,10) & (4,25,33) & (8,13,9) &  \\
\hline
(0,14,28) & (1,29,15) & (2,33,36) & \multirow{2}{*}{$\left<\alpha^{7}\right>$} \\
(3,18,41) & (6,35,26) & (9,11,38) &  \\
\hline
(0,19,26) & (1,32,17) & (2,35,25) & \multirow{5}{*}{$1$} \\
(3,24,33) & (4,7,29) & (5,41,10) &  \\
(6,8,18) & (9,14,34) & (11,31,28) &  \\
(12,30,22) & (13,23,21) & (15,40,39) &  \\
(16,27,20) & (36,37,38) &  &  \\
\hline
(0,6,12,18,24,30,36) & (1,31,19,7,37,25,13) & (2,26,8,32,14,38,20) & \multirow{2}{*}{$\left<\alpha^{3}\right>$} \\
(3,9,15,21,27,33,39) & (4,40,34,28,22,16,10) & (5,29,11,35,17,41,23) &  \\
\hline
\end{array}}\]
\end{exb}

\begin{exb}\label{hw42r39s2} $\mathrm{HW}^*(42;3^{39},7^{2})$. There are 26 orbits of length 21 on 3-cycles, and 4 orbits of length 3 on 7-cycles.
\[\fitarray{\begin{array}{|l l l |c|}
\hline
(0,10,27) & (1,14,5) &  & \multirow{6}{*}{$\left<\alpha^{3}\right>$} \\
(0,26,7) & (3,10,41) &  &  \\
(0,29,9) & (1,38,34) &  &  \\
(0,1,33) & (2,4,23) &  &  \\
(0,7,41) & (2,34,33) &  &  \\
(0,14,11) & (1,3,28) &  &  \\
\hline
(0,30,13) & (1,15,3) & (2,8,16) & \multirow{5}{*}{$1$} \\
(4,38,7) & (5,10,34) & (6,26,31) &  \\
(9,28,23) & (11,37,18) & (12,33,27) &  \\
(14,36,40) & (17,20,32) & (19,35,41) &  \\
(21,29,39) & (22,25,24) &  &  \\
\hline
(0,18,36,12,30,6,24) & (1,19,37,13,31,7,25) &  & \multirow{2}{*}{$C_{21}$} \\
(0,36,30,24,18,12,6) & (1,13,25,37,7,19,31) &  &  \\
\hline
\end{array}}\]
\end{exb}

\begin{exb}\label{hw42r40s1} $\mathrm{HW}^*(42;3^{40},7^{1})$. There are 26 orbits of length 21 and 2 orbits of length 7 on 3-cycles, and 2 orbits of length 3 on 7-cycles.
\[\fitarray{\begin{array}{|l l l |c|}
\hline
(0,28,14) & (1,15,29) &  & C_{21} \\
\hline
(0,11,3) & (1,8,40) &  & \multirow{6}{*}{$\left<\alpha^{3}\right>$} \\
(0,37,22) & (2,15,23) &  &  \\
(0,29,7) & (2,21,34) &  &  \\
(0,9,34) & (1,2,41) &  &  \\
(0,40,25) & (2,23,3) &  &  \\
(0,22,27) & (1,41,2) &  &  \\
\hline
(0,16,12) & (1,6,31) & (2,26,18) & \multirow{5}{*}{$1$} \\
(3,27,13) & (4,37,28) & (5,14,29) &  \\
(7,33,39) & (8,15,19) & (9,25,20) &  \\
(10,24,41) & (11,30,23) & (17,40,21) &  \\
(22,34,35) & (32,36,38) &  &  \\
\hline
(0,6,12,18,24,30,36) & (1,37,31,25,19,13,7) &  & C_{21} \\
\hline
\end{array}}\]
\end{exb}

\subsection*{C. Two factorizations of \texorpdfstring{$K^*_{(21:4)}$}{K*(21:4)}}

Here $x_j$ denotes the vertex $(x,j)$ of $\mathbb{Z}_{21} \times \mathbb{Z}_4$. In each example, the first table lists the rows of base $3$-cycles, each developed by the subgroup $\left<\alpha^3\right>$ of order $7$ and yielding three $\vv{C}_3$-factors, and the second table lists the rows of base $7$-cycles, each developed by the identity subgroup and yielding $21$ $\vv{C}_7$-factors.

\begin{exc}\label{k214a21} A $\{\vv{C}_3^{\,21}, \vv{C}_7^{\,42}\}$-factorization of $K^*_{(21:4)}$ on $\mathbb{Z}_{21} \times \mathbb{Z}_4$ developed by $C_{21}=\left<\alpha\right>$ where $\alpha~:~x_j\rightarrow (x+1)_j$; the parts are $\mathbb{Z}_{21} \times \{j\}$, $j \in \mathbb{Z}_4$. There are 28 orbits of length 21 on 3-cycles, and 24 orbits of length 21 on 7-cycles.
\[\fitarray{\begin{array}{|l l l l|c|}
\hline
(1_{0},18_{2},20_{1}) & (0_{0},16_{3},7_{1}) & (2_{0},13_{2},5_{3}) & (0_{1},12_{3},17_{2}) & \multirow{7}{*}{$\left<\alpha^{3}\right>$} \\
(0_{0},13_{1},14_{2}) & (2_{0},19_{3},14_{1}) & (1_{0},10_{2},3_{3}) & (0_{1},5_{3},12_{2}) &  \\
(1_{0},12_{1},16_{2}) & (0_{0},8_{1},7_{3}) & (2_{0},3_{3},2_{2}) & (1_{1},17_{3},15_{2}) &  \\
(1_{0},19_{2},16_{1}) & (0_{0},4_{3},2_{1}) & (2_{0},8_{2},11_{3}) & (0_{1},3_{3},6_{2}) &  \\
(1_{0},13_{1},11_{2}) & (2_{0},2_{3},11_{1}) & (0_{0},18_{3},19_{2}) & (0_{1},1_{3},9_{2}) &  \\
(2_{0},8_{1},13_{2}) & (1_{0},15_{3},15_{1}) & (0_{0},12_{2},13_{3}) & (1_{1},17_{2},8_{3}) &  \\
(1_{0},17_{1},8_{2}) & (2_{0},16_{1},20_{3}) & (0_{0},7_{2},4_{3}) & (0_{1},0_{2},6_{3}) &  \\
\hline
\end{array}}\]
\setlength{\arraycolsep}{3pt}
{\small\[\fitarray{\begin{array}{|l l l|c|}
\hline
(0_{0},4_{2},0_{3},7_{0},11_{1},1_{2},5_{3}) & (1_{0},8_{1},3_{0},11_{3},17_{0},12_{2},4_{1}) & (2_{0},17_{3},7_{2},12_{1},11_{2},1_{1},18_{2}) & \multirow{4}{*}{$1$} \\
(4_{0},6_{3},14_{1},3_{2},19_{1},11_{0},9_{2}) & (5_{0},2_{1},19_{0},1_{3},3_{1},9_{3},6_{2}) & (6_{0},19_{3},2_{2},0_{1},15_{0},14_{2},4_{3}) &  \\
(8_{0},20_{3},10_{1},19_{2},16_{0},0_{2},7_{3}) & (9_{0},7_{1},18_{3},20_{2},5_{1},13_{0},9_{1}) & (10_{0},15_{3},10_{2},6_{1},16_{3},16_{2},10_{3}) &  \\
(12_{0},13_{1},5_{2},20_{0},13_{2},14_{0},8_{2}) & (18_{0},17_{1},3_{3},15_{2},13_{3},16_{1},8_{3}) & (15_{1},2_{3},20_{1},17_{2},12_{3},18_{1},14_{3}) &  \\
\hline
(0_{0},11_{3},1_{0},6_{1},9_{0},17_{2},16_{3}) & (2_{0},8_{3},14_{2},12_{0},1_{1},11_{0},10_{3}) & (3_{0},12_{1},19_{2},6_{0},4_{3},11_{1},9_{3}) & \multirow{4}{*}{$1$} \\
(4_{0},14_{3},13_{1},6_{3},15_{0},3_{3},20_{1}) & (5_{0},5_{2},7_{3},1_{2},10_{0},4_{1},1_{3}) & (7_{0},9_{2},17_{3},16_{0},19_{1},19_{3},15_{2}) &  \\
(8_{0},8_{1},10_{2},18_{1},17_{0},18_{2},7_{1}) & (13_{0},20_{3},9_{1},12_{2},5_{1},20_{2},2_{1}) & (14_{0},6_{2},15_{3},4_{2},3_{1},18_{3},11_{2}) &  \\
(18_{0},7_{2},12_{3},0_{2},0_{3},13_{2},14_{1}) & (19_{0},0_{1},2_{3},15_{1},8_{2},20_{0},2_{2}) & (10_{1},16_{2},16_{1},3_{2},13_{3},17_{1},5_{3}) &  \\
\hline
\end{array}}\]}
\end{exc}

\begin{exc}\label{k214a42} A $\{\vv{C}_3^{\,42}, \vv{C}_7^{\,21}\}$-factorization of $K^*_{(21:4)}$ on $\mathbb{Z}_{21} \times \mathbb{Z}_4$ developed by $C_{21}=\left<\alpha\right>$ where $\alpha~:~x_j\rightarrow (x+1)_j$; the parts are $\mathbb{Z}_{21} \times \{j\}$, $j \in \mathbb{Z}_4$. There are 56 orbits of length 21 on 3-cycles, and 12 orbits of length 21 on 7-cycles.
\[\fitarray{\begin{array}{|l l l l|c|}
\hline
(0_{0},3_{1},8_{2}) & (2_{0},14_{1},4_{3}) & (1_{0},6_{2},20_{3}) & (1_{1},7_{2},18_{3}) & \multirow{14}{*}{$\left<\alpha^{3}\right>$} \\
(0_{0},2_{2},5_{1}) & (1_{0},5_{3},13_{1}) & (2_{0},1_{2},16_{3}) & (0_{1},0_{3},9_{2}) &  \\
(1_{0},17_{1},13_{2}) & (0_{0},19_{1},3_{3}) & (2_{0},6_{2},19_{3}) & (0_{1},8_{3},11_{2}) &  \\
(1_{0},19_{1},6_{2}) & (2_{0},18_{3},12_{1}) & (0_{0},8_{2},20_{3}) & (2_{1},16_{2},13_{3}) &  \\
(1_{0},6_{1},7_{2}) & (2_{0},13_{1},14_{3}) & (0_{0},12_{2},15_{3}) & (2_{1},1_{3},17_{2}) &  \\
(2_{0},10_{1},9_{2}) & (1_{0},20_{3},0_{1}) & (0_{0},1_{2},1_{3}) & (2_{1},20_{2},18_{3}) &  \\
(1_{0},8_{2},3_{1}) & (0_{0},13_{1},16_{3}) & (2_{0},16_{2},2_{3}) & (2_{1},15_{3},12_{2}) &  \\
(1_{0},3_{1},15_{2}) & (0_{0},13_{3},11_{1}) & (2_{0},11_{3},13_{2}) & (1_{1},18_{3},17_{2}) &  \\
(2_{0},12_{2},10_{1}) & (1_{0},15_{1},12_{3}) & (0_{0},2_{3},13_{2}) & (2_{1},11_{2},16_{3}) &  \\
(2_{0},19_{2},6_{1}) & (1_{0},7_{3},16_{1}) & (0_{0},12_{3},20_{2}) & (2_{1},6_{2},2_{3}) &  \\
(2_{0},18_{2},18_{1}) & (0_{0},14_{3},13_{1}) & (1_{0},1_{3},20_{2}) & (2_{1},9_{3},19_{2}) &  \\
(0_{0},7_{1},18_{2}) & (2_{0},6_{1},10_{3}) & (1_{0},18_{3},10_{2}) & (2_{1},11_{3},5_{2}) &  \\
(1_{0},16_{2},8_{1}) & (0_{0},15_{3},18_{1}) & (2_{0},5_{3},12_{2}) & (1_{1},11_{2},10_{3}) &  \\
(1_{0},14_{2},15_{1}) & (0_{0},20_{1},18_{3}) & (2_{0},7_{3},3_{2}) & (1_{1},16_{2},5_{3}) &  \\
\hline
\end{array}}\]
\setlength{\arraycolsep}{3pt}
{\small\[\fitarray{\begin{array}{|l l l|c|}
\hline
(0_{0},11_{3},1_{0},7_{2},3_{0},12_{2},0_{1}) & (2_{0},2_{1},17_{3},2_{2},8_{3},3_{1},9_{3}) & (4_{0},5_{3},16_{1},20_{0},5_{1},8_{2},10_{3}) & \multirow{4}{*}{$1$} \\
(5_{0},16_{2},20_{3},0_{2},19_{0},19_{2},14_{3}) & (6_{0},4_{2},19_{1},12_{3},8_{0},5_{2},12_{1}) & (7_{0},10_{2},15_{0},4_{1},17_{2},18_{3},1_{2}) &  \\
(9_{0},18_{1},20_{2},7_{3},9_{1},4_{3},9_{2}) & (10_{0},11_{1},11_{2},10_{1},1_{3},15_{1},13_{2}) & (11_{0},7_{1},14_{2},18_{0},15_{3},15_{2},3_{3}) &  \\
(12_{0},6_{1},16_{3},1_{1},13_{0},0_{3},13_{1}) & (14_{0},13_{3},6_{2},8_{1},3_{2},17_{1},19_{3}) & (16_{0},2_{3},20_{1},17_{0},6_{3},18_{2},14_{1}) &  \\
\hline
\end{array}}\]}
\end{exc}


\begin{thebibliography}{99}

\bibitem{Abel2002} R.J.R. Abel, F.E. Bennett and G. Ge, {\it Resolvable perfect Mendelsohn designs with block size five}, Discrete Math. \textbf{247}, 1--12, 2002.

\bibitem{Adamsun} P. Adams and D. Bryant, {\it Resolvable directed cycle systems of all indices for cycle length 3 and 4}, unpublished.

\bibitem{ASSW} B. Alspach, P.J. Schellenberg, D.R. Stinson and D. Wagner, {\it The Oberwolfach problem and factors of uniform odd length cycles}, J. Combin. Theory Ser. A \textbf{52}, 20--43, 1989.

\bibitem{AssafHartman} A. Assaf and A. Hartman, {\it Resolvable group divisible designs with block size 3}, Discrete Math. \textbf{77}, 5--20, 1989.

\bibitem{Bennett1990} F.E. Bennett and X. Zhang, {\it Resolvable Mendelsohn designs with block size $4$}, Aequationes Math. \textbf{40}, 248--260, 1990.

\bibitem{Bermond1979} J.C. Bermond, A. Germa and D. Sotteau, {\it Resolvable decomposition of $K_{n}^{*}$}, J. Combin. Theory Ser. A \textbf{26}, 179--185, 1979.

\bibitem{Sajna2018} A. Burgess, N. Franceti\'c and M. \v{S}ajna, {\it On the directed Oberwolfach problem with equal cycle lengths: the odd case}, Australas. J. Combin. \textbf{71}, 272--292, 2018.

\bibitem{Sajna2014} A. Burgess and M. \v{S}ajna, {\it On the directed Oberwolfach problem with equal cycle lengths}, Electron. J. Combin. \textbf{21}, P1.15, 2014.

\bibitem{DPR} I.J. Dejter, D.A. Pike and C.A. Rodger, {\it The directed almost resolvable Hamilton--Waterloo problem}, Australas. J. Combin. \textbf{18}, 201--208, 1998.

\bibitem{Francetic2023} N. Franceti\'c and M. \v{S}ajna, {\it On the directed Oberwolfach problem for complete symmetric equipartite digraphs and uniform-length cycles}, J. Combin. Des. \textbf{31}, 604--641, 2023.

\bibitem{HoffmanSchellenberg} D.G. Hoffman and P.J. Schellenberg, {\it The existence of $C_k$-factorizations of $K_{2n}-F$}, Discrete Math. \textbf{97}, 243--250, 1991.

\bibitem{Lacaze2023} A. Lacaze-Masmonteil, {\it Completing the solution of the directed Oberwolfach problem with cycles of equal length}, J. Combin. Des. \textbf{32}, 5--30, 2024.

\bibitem{LeiFu} H.C. Lei and H.L. Fu, {\it The Hamilton--Waterloo problem for triangle-factors and heptagon-factors}, Graphs Combin. \textbf{32}, 271--278, 2016.

\bibitem{Liu2000} J. Liu, {\it A generalization of the Oberwolfach problem and $C_t$-factorizations of complete equipartite graphs}, J. Combin. Des. \textbf{8}, 42--49, 2000.

\bibitem{Liu2003} J. Liu, {\it The equipartite Oberwolfach problem with uniform tables}, J. Combin. Theory Ser. A \textbf{101}, 20--34, 2003.

\bibitem{Piotrowski} W.L. Piotrowski, {\it The solution of the bipartite analogue of the Oberwolfach problem}, Discrete Math. \textbf{97}, 339--356, 1991.

\bibitem{RCW} D.K. Ray-Chaudhuri and R.M. Wilson, {\it Solution of Kirkman's schoolgirl problem}, in: Combinatorics (Proc. Sympos. Pure Math., Vol. XIX), 187--203, Amer. Math. Soc., Providence, RI, 1971.

\bibitem{Rees93} R. Rees, {\it Two new direct product-type constructions for resolvable group-divisible designs}, J. Combin. Des. \textbf{1}, 15--26, 1993.

\bibitem{Tillson1980} T.W. Tillson, {\it A Hamiltonian decomposition of $K_{2m}^*$, $2m \geq 8$}, J. Combin. Theory Ser. B \textbf{29}, 68--74, 1980.

\bibitem{WCC} L. Wang, F. Chen and H. Cao, {\it The Hamilton--Waterloo problem for $C_3$-factors and $C_n$-factors}, J. Combin. Des. \textbf{25}, 385--418, 2017.

\bibitem{YOK} F. Yetgin, U. Odaba\c{s}{\i} and E. Koloto\u{g}lu, {\it The directed Hamilton--Waterloo problem with $4$-cycles and odd cycles}, submitted.

\bibitem{FUS} F. Yetgin, U. Odaba\c{s}{\i} and S. \"Ozkan, {\it On the directed Hamilton--Waterloo problem with two cycle sizes}, Contrib. Discrete Math. \textbf{20}, 74--94, 2025.

\bibitem{FUS2} F. Yetgin, U. Odaba\c{s}{\i} and S. \"Ozkan, {\it The directed uniform Hamilton--Waterloo problem involving even cycle sizes}, Discuss. Math. Graph Theory \textbf{45}, 615--636, 2025.

\end{thebibliography}
\end{document}